\documentclass[lettersize,journal]{IEEEtran}
\usepackage{amsmath,amsfonts}
\usepackage{algorithmic}
\usepackage{algorithm}
\usepackage{array}
\usepackage[caption=false,font=normalsize,labelfont=sf,textfont=sf]{subfig}
\usepackage{textcomp}
\usepackage{stfloats}
\usepackage{url}
\usepackage{verbatim}
\usepackage{graphicx}
\usepackage{cite}
\usepackage{colortbl}
\usepackage{amssymb}

\newtheorem{theorem}{Theorem}

\newtheorem{remark}{Remark}
\newtheorem{definition}{Definition}
\newtheorem{assumption}{Assumption}
\newtheorem{lemma}{Lemma}
\newtheorem{example}{Example}
\begin{document}

\title{Weak Closed-loop Solvability for Discrete-time Stochastic Linear-Quadratic Optimal Control}

\author{Yue Sun, Xianping Wu and Xun Li ~%\IEEEmembership{***,~***,}
        % <-this % stops a space
\thanks{This research is supported by the National Natural Science Foundation of China (No. 12201129)
and Natural Science Foundation of Guangdong Province
(No. 2022A1515010839).

 }% <-this % stops a space
\thanks{ Y. Sun is with the School of Control Science and Engineering, Shandong University, Jinan, Shandong 250061, China . X. Wu is with the School of  Mathematics and Statistics, Guangdong University of Technology, Guangzhou, China (xianpingwu@gdut.edu.cn). X. Li is with the Department of Applied Mathematics, The Hong Kong Polytechnic University, Hong Kong, China}}



\maketitle

\begin{abstract}
In this paper, the solvability of discrete-time stochastic linear-quadratic (LQ) optimal control problem in finite horizon is considered.
Firstly, it shows that the closed-loop solvability for the LQ control problem is optimal if and only if the generalized Riccati equation admits a regular solution by solving the forward and backward difference equations iteratively.
To this ends, it finds that the open-loop solvability is strictly weaker than closed-loop solvability, that is, the LQ control problem is always open-loop optimal solvable if it is closed-loop optimal solvable but not vice versa.
Secondly, by the perturbation method, it proves that the weak-closed loop strategy which is a feedback form of a state feedback representation is equivalent to the open-loop solvability of the LQ control problem.
Finally, an example sheds light on the theoretical results established.
\end{abstract}


\section{Introduction}

Optimal control theory is an important component of modern control theory. As early as 1948, Wiener \cite{Wiener1948} introduced concepts such as control and feedback, laying an important foundation for control theory.
Pioneered by Kalman \cite{Kalman1960} in 1960, linear-quadratic (LQ) optimal control has been playing a classical role in modern control theory.
However, considering the uncertainty and inevitable random interference factors in reality, Wonham \cite{Wonham1968} first extended LQ control theory to stochastic systems in 1968.
In 1976, Bismut \cite{Bismut1976} solved a stochastic system with multiplicative noise and discussed the standard stochastic LQ optimal control problem with random coefficients using functional analysis techniques. Based on the Riccati equation, the existence of the optimal controller was obtained.
From then on, many scholars began to study stochastic optimal control problems thereafter. A common assumption in most literature on stochastic LQ problems is that the state weighting matrix is non negative definite and the control weighting matrix is positive definite.
That is to say, consider the following system
\begin{eqnarray*}
  &\hspace{-0.8em}&\hspace{-0.8em}\min \sum^{N-1}_{t=0}\mathbb{E}(x'_tQx_t+2u'_tSx_t+u'_tRu_t)+\mathbb{E}(x'_{N}Hx_{N}),\\
  &\hspace{-0.8em}&\hspace{-0.8em}\mbox{subject\!\quad\!to}\quad  x_{t+1}=(A+w_tC)x_t+(B+w_tD)u_t,
\end{eqnarray*}
where $x_t\in \mathbb{R}^n$ is the state with the initial value $x_0$, $u_t\in \mathbb{R}^m$ is the control input, $A$, $B$, $C$  and $D$ are constant matrices with compatible dimensions.
Under the traditional assumption, the cost functional always has positive semi-definite weighting matrices for the control and the state. For a constant $\delta >0$,
%assume that  $\forall t\in{\color{blue} \mathbb{N}=\{0, 1, \cdots, N-1\}}$,
assume that
\begin{eqnarray*}
  R\geq \delta I, \quad Q-S'R^{-1}S\geq 0.
\end{eqnarray*}
Under these conditions, the LQ problem admits a unique solution if and only if the Riccati equation with the terminal value $P_{N}=H$ has a unique positive definite solution satisfying
\begin{eqnarray}\label{0}
P_t&\hspace{-0.8em}=&\hspace{-0.8em}Q+A'P_{t+1}A+C'P_{t+1}C-(A'P_{t+1}B+C'P_{t+1}D\nonumber\\
&\hspace{-0.8em}&\hspace{-0.8em}+S')
(R+B'P_{t+1}B+D'P_{t+1}D)^{-1}\nonumber\\
&\hspace{-0.8em}&\hspace{-0.8em}\times (B'P_{t+1}A+D'P_{t+1}C+S),
\end{eqnarray}
while the unique open-loop optimal solution to the problem is a linear feedback form of the state.
Recently, Zhang \cite{Zhang2015} studied a class of discrete-time stochastic systems involving multiplicative noises and input delay where the optimal control was shown to be a linear function of the conditional expectation of the state, with the feedback gain based on a Riccati-ZXL difference equation. Qi in \cite{Qi2018} and \cite{Qi2022} successively investigated the mean-field stochastic system without or with input delay by convexity condition, while the optimal control was derived in terms of the solution of a modified Riccati equation.

Contrary to the traditional positive-definite assumption,  Zhou revealed a fact that even if the cost weighting matrices are indefinite the stochastic LQ optimal control problem with multiplicative noises might still be solvable in \cite{Zhou1988}.
Moreover, the general necessary and sufficient conditions for the solvability of the generalized differential Riccati equation with an indefinite quadratic cost function for a linear It$\hat{o}$ system subject to multiplicative noise on both the state and control was derived in \cite{Zhou2001}. For discrete-time system, a  Riccati equation was obtained  to characterize the existence and the form of the optimal control:
\begin{eqnarray}\label{7}
P_t
&\hspace{-0.8em}=&\hspace{-0.8em} Q+A'P_{t+1}A+C'P_{t+1}C-(A'P_{t+1}B+C'P_{t+1}D\nonumber\\
&\hspace{-0.8em}&\hspace{-0.8em}+S')({R}+B'P_{t+1}B+D'P_{t+1}D)^{\dagger}\nonumber\\
&\hspace{-0.8em}&\hspace{-0.8em}\times (B'P_{t+1}A+D'P_{t+1}C+S),
\end{eqnarray}
with the terminal value $P_N=H$ where it didn't not require any definiteness constraint on the cost weighting matrices.
Inspired by \cite{Sun2014} and \cite{Sun2016}, which focused on the open-loop and closed-loop solvability of the stochastic It$\hat{o}$ system, the equivalent relationship between the open-loop and closed-loop solvability is considered. As is known to all, under the assumption of positive definite condition of weighted matrix, they are equivalent. For the indefinite LQ problem, it has found that the closed-loop solvability implies the open-loop solvability, however, there are some LQ problems which is open-loop solvable rather than closed-loop solvable, that is the optimal solution cannot be obtained from the Riccati equation.
More details can be seen from the following example:
\begin{example}\label{example1}
Consider the following system with one-dimension state equation:
\begin{eqnarray}\label{01}
  x_{t+1} =x_t+u_t+w_tu_t, \quad t\in\{0,1\},
\end{eqnarray}
with the initial value $x_0$. $w_t$ is the Gaussian noise with $\mathbb{E}[w_t]=0$ and $\mathbb{E}[w_tw_l]=\delta_{kl}$. $x_0$ and $w_t$ are independent with each other.  The associated cost function is given by
\begin{eqnarray}\label{02}
  J(x_0, u) &\hspace{-0.8em}=&\hspace{-0.8em} \mathbb{E}\Big(x^2_2-2\sum\limits^{1}_{t=0}u^2_t\Big).
\end{eqnarray}

Following from the difference Riccati equation
\begin{eqnarray*}
P_t&\hspace{-0.8em}=&\hspace{-0.8em} Q+A'P_{t+1}A+C'P_{t+1}C-(A'P_{t+1}B+C'P_{t+1}D) \nonumber\\
&\hspace{-0.8em}&\hspace{-0.8em}\times(R+B'P_{t+1}B+D'P_{t+1}D)^{\dagger}(B'P_{t+1}A+D'P_{t+1}C),
\end{eqnarray*}
with the terminal value $P_2=1$, and by simple calculation, the solution of it is exactly $P_t\equiv 1$, $t\in \{0, 1, 2\}$.

Combining with the maximum principle, a usual Riccati equation approach specifies the corresponding state feedback control as follows
\begin{eqnarray}\label{03}
  u^*_t
&\hspace{-0.8em}=&\hspace{-0.8em} -(R+B'P_{t+1}B+D'P_{t+1}D)^{\dagger}\nonumber\\
&\hspace{-0.8em}&\hspace{-0.8em}\times (B'P_{t+1}Ax_t+D'P_{t+1}C) \equiv0, \quad t\in\{0, 1\},
\end{eqnarray}
due to $R=-2$ and $0^{\dagger}=0$.
It should be noted that (\ref{03}) is no longer an optimal open-loop solution if $x_0\neq 0$. To be specific, adding $u^*_t$ in (\ref{03}) into (\ref{01}), the state $x^*_2$ is calculated as
\begin{eqnarray*}
  x^*_2 &\hspace{-0.8em}=&\hspace{-0.8em} x_0 +B_0u^*_0+B_1u^*_1=x_0,
\end{eqnarray*}
where $B_t=1+w_t$.

To this end, the cost function in (\ref{02}) is given by
\begin{eqnarray*}
  J(x_0, u^*) &\hspace{-0.8em}=&\hspace{-0.8em}\mathbb{E}(x^2_0) >0,
\end{eqnarray*}
for any $x_0\neq 0$.

Moreover, let $\bar{u}_t$ be the control defined by
\begin{eqnarray*}
 \bar{u}_t=Sx_t, \quad S=-0.3194.
\end{eqnarray*}
In this way, we have
\begin{eqnarray*}
\bar{u}_0&\hspace{-0.8em}=&\hspace{-0.8em}Sx_0,\\
 \bar{u}_1&\hspace{-0.8em}=&\hspace{-0.8em}Sx_1=S(x_0+B_0u_0)=S(x_0+SB_0x_0)\nonumber\\
&\hspace{-0.8em}=&\hspace{-0.8em}S(1+SB_0)x_0.
\end{eqnarray*}
Adding them into (\ref{01}), $x_2$ is rewritten as
\begin{eqnarray*}
  x_2 &\hspace{-0.8em}=&\hspace{-0.8em} x_0+B_0\bar{u}_0+B_1\bar{u}_1.
\end{eqnarray*}
Thus, the cost function $J(x_0, \bar{u})$ is given by
\begin{eqnarray}\label{010}
J(x_0, \bar{u})&\hspace{-0.8em}=&\hspace{-0.8em}\mathbb{E}(x^2_2-2\bar{u}^2_0-2\bar{u}^2_1)\nonumber\\
&\hspace{-0.8em}=&\hspace{-0.8em}\mathbb{E}(x_0+B_0\bar{u}_0+B_1\bar{u}_1)^2-2E(\bar{u}^2_0)-2E(\bar{u}^2_1)\nonumber\\
&\hspace{-0.8em}=&\hspace{-0.8em}\mathbb{E}(x^2_0+2x_0B_0\bar{u}_0+2x_0B_1\bar{u}_1
+\bar{u}_0B_0B_0\bar{u}_0\nonumber\\
&\hspace{-0.8em}&\hspace{-0.8em}+2\bar{u}_0B_0B_1\bar{u}_1
+\bar{u}_1B_1B_1\bar{u}_1)-2\mathbb{E}(\bar{u}^2_0)-2\mathbb{E}(\bar{u}^2_1)\nonumber\\
&\hspace{-0.8em}=&\hspace{-0.8em}(4S^3+4S^2+4S+1)\mathbb{E}(x^2_0),%\nonumber\\
%&\hspace{-0.8em}=&\hspace{-0.8em}0<J(x_0, u^*_t),
\end{eqnarray}
where
\begin{eqnarray*}
  \mathbb{E}(\bar{u}_iB_iB_i\bar{u}_i)&\hspace{-0.8em}=&\hspace{-0.8em}2\mathbb{E}(\bar{u}^2_i), \quad i=0, 1,\\
  \mathbb{E}(x_0B_0\bar{u}_0)&\hspace{-0.8em}=&\hspace{-0.8em}S\mathbb{E}(x^2_0),\\
  \mathbb{E}(x_0B_1\bar{u}_1)&\hspace{-0.8em}=&\hspace{-0.8em}\mathbb{E}[x_0B_1S(1+SB_0)x_0]
  =S\mathbb{E}(x^2_0)\nonumber\\
  &\hspace{-0.8em}&\hspace{-0.8em}+S^2\mathbb{E}(x^2_0),\\
  \mathbb{E}(\bar{u}_0B_0B_1\bar{u}_1)&\hspace{-0.8em}=&\hspace{-0.8em}\mathbb{E}
  [x_0SB_0B_1S(1+SB_0)x_0]\nonumber\\
  &\hspace{-0.8em}=&\hspace{-0.8em}\mathbb{E}(x_0S^2B_0B_1x_0+x_0S^3B_0B_1B_0x_0)\\
  &\hspace{-0.8em}=&\hspace{-0.8em}S^2\mathbb{E}(x^2_0)+2S^3E(x^2_0)
\end{eqnarray*}
have been used in the derivation.

Since the real root of $4S^3+4S^2+4S+1=0$ is calculated as
\begin{eqnarray}
  S&\hspace{-0.8em}=&\hspace{-0.8em}\sqrt[3]{\frac{1}{216}+\sqrt{\frac{1}{216^2}+\frac{2^3}{9^3}}}
  + \sqrt[3]{\frac{1}{216}-\sqrt{\frac{1}{216^2}+\frac{2^3}{9^3}}}-\frac{1}{3}\nonumber\\
  &\hspace{-0.8em}=&\hspace{-0.8em}-0.3194,
\end{eqnarray}
by letting $\bar{u}_t=Sx_t$ and $S=-0.3194$, for any $x_0 \neq 0$, (\ref{010}) degrades into the following form
\begin{eqnarray*}
J(x_0, \bar{u})&\hspace{-0.8em}=&\hspace{-0.8em}0\cdot \mathbb{E}(x^2_0)\nonumber\\
&\hspace{-0.8em}=&\hspace{-0.8em}0<J(x_0, u^*).
\end{eqnarray*}

Based on the discussion above, since the cost function is nonnegative, $\bar{u}_t$ is open-loop optimal for the initial value $x_0$, but $u^*_t$ is not.
\end{example}

Following from the Example, it concludes that there exist some stochastic LQ control problems where the usual solvability of the generalized Riccati equation  may not be helpful in dealing with the open-loop solvability. Naturally, there generates a question that when the stochastic LQ problem is only open-loop solvable, does it still exist a  feedback form of the state for an open-loop optimal control? Out of consideration for this point and motivated by \cite{Sun2018} and \cite{Wen2021}, it presented a weak closed-loop solution whose result is exactly an open-loop solvability for stochastic It$\hat{o}$ differential system, in this paper, firstly, the equivalent relationship for the closed-loop solvability of the stochastic discrete-time system without  any definiteness constraint is given.
Secondly, using the perturbation approach proposed in \cite{Sun2016}, the equivalences of open-loop solvability are derived which lays the foundation for the discussion the weak closed-loop solvability satisfying a linear feedback form of the sate.
Finally, we will shown that the stochastic LQ problem is open-loop solvable if and only if there always exists a weak closed-loop strategy whose outcome is an open-loop optimal control.

The rest of the paper is organized as follows. Problem formulation and the closed-loop solvability related with the Riccati equation is arranged in Section 2.
The open-loop solvability by perturbation method are presented in Section 3.
In Section 4, we are devoted to show the process to obtain a weak closed-loop solvability and prove the equivalent relationship between the open-loop solvability and the weak closed-loop solvability. A numerical example in Section 5 sheds light on the results established in
this paper.

\emph{Notations:} $\mathbb{R}^{m\times 1}$ is the sets of $m$-dimensional Euclidean space.
$\{\Omega, \mathcal{F}, \mathbb{P}, \{\mathcal{F}_t\}_{t \geq 0}\}$  is a complete probability space on which a scalar white noise $w_t $ is defined such that
$\{\mathcal{F}_t\}_{t\geq 0}$  is the natural filtration generated by $w_t $, i.e., $\mathcal{F}_t=\sigma\{w_0, ... ,w_t\}$.
The conditional mathematical expectation with regard to the filtration $\mathcal{F}_t$ is denoted as $\mathbb{E}[\cdot| \mathcal{F}_t]$.
$A'$ is the transpose of the matrix $A$.
A symmetric matrix $A>0$ (or $\geq0$) indicates that $A$ is a positive definite matrix (or positive semi-definite matrix).
$M^{\dag}$ represents the Moore-Penrose pseudoinverse of a matrix $M$ satisfying $MM^{\dag}M=M$, $M^{\dag}MM^{\dag}=M^{\dag}$, $(MM^{\dag})'=MM^{\dag}$ and $(M^{\dag}M)'=M^{\dag}M$.
$Rang(A)=\{Ax| x\in \mathbb{R}^n\}$ with $n$ being the dimension of $x$ is the range of $A$.
Let $\mathbb{N}$ denote the set of $\{0, 1, \cdots, N\}$.
$L^2_{\mathcal{F}}(a; \mathbb{R}^n)$ denotes the space of all $\mathbb{R}^n$-valued $\mathcal{F}_{a-1}$-adapted  variables satisfying $\{\varphi: \Omega\rightarrow \mathbb{R}^n| \mathbb{E}|\varphi|^2 < \infty\}$. $L^2_{\mathcal{F}}(\mathbb{N}; \mathbb{R}^m)$ denotes the space of all $\mathbb{R}^m$-valued $\mathcal{F}_{t-1}$-adapted  variables satisfying $\{u: \mathbb{N}\times \Omega\rightarrow \mathbb{R}^m\Big| \mathbb{E}\sum\limits^{N-1}_{t=0} |u_t|^2 < \infty\}$.
$L^2(\mathbb{N}; \mathbb{R}^m)$ denotes the space of all $\mathbb{R}^m$-valued deterministic variables satisfying $\{\phi: \mathbb{N} \rightarrow \mathbb{R}^m\Big| \sum\limits^{N-1}_{t=0} |\phi_t|^2 < \infty\}$.
%$L^2(0, N; \mathbb{R}^m)=\{[0, N] \rightarrow \mathbb{R}^m\Big|  \sum\limits^{N-1}_{t=0} |u_t|^2 < \infty\}$

\section{Problem Formulation}
Consider the discrete-time system:
\begin{eqnarray}\label{1}
  x_{t+1} &\hspace{-0.8em}=&\hspace{-0.8em} A_tx_t+B_tu_t+b_t+(C_tx_t+D_tu_t+\sigma_t)w_t,
\end{eqnarray}
where $x_t\in \mathbb{R}^n$ is the state with the initial value $x_0$. $u_t\in \mathbb{R}^m$ is the control input. $A_t$, $B_t$ , $C_t$ and $D_t$ are deterministic matrix-valued functions of compatible dimensions. $b_t$ and $\sigma_t$ are $\mathcal{F}_{t-1}$-adapted.
%, where $\mathcal{F}_t=\{w_0, ..., w_t\}$.
The associated cost function is given by
\begin{eqnarray}\label{2}
  J(x_0, u) &\hspace{-0.8em}=&\hspace{-0.8em} \mathbb{E}\Big(\sum\limits^{N-1}_{t=0}[x'_tQ_tx_t+2u'_tS_tx_t+u'_tR_tu_t+2x'_tq_t\nonumber\\
  &\hspace{-0.8em}&\hspace{-0.8em}+2u'_t \rho_t]+x'_{N}Hx_{N}+2x'_{N}g\Big),
\end{eqnarray}
where $H\in \mathbb{R}^{n\times n}$ is a symmetric constant matrix. $Q_t\in \mathbb{R}^{n\times n}$, $R_t\in \mathbb{R}^{m\times m}$ and $S_t\in \mathbb{R}^{m\times n}$ are deterministic matrix-valued functions with $Q_t$ and $R_t$ being symmetric. $g$ is $\mathcal{F}_{N-1}$-adapted random variable. $q_t$ and $\rho_t$ are $\mathcal{F}_{t-1}$-adapted.

The problem studied in this paper is given.
\smallskip

\emph{Problem (SLQ).} For any given initial value $x_0\in \mathbb{R}^n$, find a controller
$\bar{u}\in L^2_{\mathcal{F}}(\mathbb{N}; \mathbb{R}^{m})$ satisfying
\begin{eqnarray*}
J(x_0, \bar{u}) \leq J(x_0, u), \quad \forall u\in L^2_{\mathcal{F}}(\mathbb{N}; \mathbb{R}^{m}).
\end{eqnarray*}
%where $\mathcal{U}[0, N]$ is the admissible control set which is denoted as
%\begin{eqnarray*}
%% \nonumber to remove numbering (before each equation)
% \mathcal{U}[0, N]&\hspace{-0.8em}=&\hspace{-0.8em}L^2_{\mathcal{F}}(0, N; \mathbb{R}^{m})\nonumber\\
%  &\hspace{-0.8em}=&\hspace{-0.8em}\{u: [0, N]\times \Omega\rightarrow \mathbb{R}^{m}\Big| E \sum\limits^{N-1}_{t=0}|u_t|^2< \infty\}.
%\end{eqnarray*}

Moreover, the value function of Problem (SLQ) with the initial value $x_0$ is defined as
\begin{eqnarray*}
  V(x_0) =\inf_{u_t\in L^2_{\mathcal{F}}(\mathbb{N}; \mathbb{R}^{m})}J(x_0, u).
\end{eqnarray*}

Noting that for the special case, i.e., when $b_t=\sigma_t=q_t=\rho_t=g=0$, (\ref{1})-(\ref{2}) degenerate into
\begin{eqnarray}
  x_{t+1} &\hspace{-0.8em}=&\hspace{-0.8em} A_tx_t+B_tu_t+(C_tx_t+D_tu_t)w_t,\label{2-1}\\
J^0(x_0, u) &\hspace{-0.8em}=&\hspace{-0.8em} \mathbb{E}\Big(\sum\limits^{N-1}_{t=0}[x'_tQ_tx_t+2u'_tS_tx_t+u'_tR_tu_t]+x'_{N}Hx_{N}\Big),\label{2-2}
\end{eqnarray}
respectively, with the initial value $x_0$. The problem of minimizing (\ref{2-2}) subject to (\ref{2-1}) is said to be a homogeneous LQ problem related with Problem (SLQ), which is called as Problem (SLQ)$^0$. And the value function of  Problem (SLQ)$^0$ is denoted as $V^0(x_0)$.

Before dealing with the problem, there are assumptions about the coefficients in (\ref{1}) and weighting coefficients in (\ref{2}).
\begin{assumption}\label{A1}
The coefficients $A: \mathbb{N}\rightarrow \mathbb{R}^{n\times n}$ and
 $B: \mathbb{N}\rightarrow \mathbb{R}^{n\times m}$ are uniformly bounded mappings.
$b$ and $\sigma$ belong to $L^2_{\mathcal{F}}(\mathbb{N}; \mathbb{R}^n)$.
\end{assumption}

\begin{assumption}\label{A2}
The coefficients $Q: \mathbb{N}\rightarrow \mathbb{R}^{n\times n}$,
$R: \mathbb{N}\rightarrow \mathbb{R}^{m\times m}$ and
$S: \mathbb{N}\rightarrow \mathbb{R}^{m\times n}$ are uniformly bounded mappings.
$q\in L^2_{\mathcal{F}}(\mathbb{N}; \mathbb{R}^n)$, $\rho \in L^2_{\mathcal{F}}(\mathbb{N}; \mathbb{R}^m)$ and $g\in L^2_{\mathcal{F}}(N-1; \mathbb{R}^n)$.
\end{assumption}

Before proposing the main results of this paper, the definitions of open-loop, closed-loop and weak closed-loop solvable of Problem (SLQ) are given.

\begin{definition}
Problem (SLQ) is (uniquely) open-loop solvable with the initial value $x_0$ if there exists a (unique) $\bar{u}=\bar{u}(x_0)\in L^2_{\mathcal{F}}(\mathbb{N}; \mathbb{R}^m)$, which means it's related with the initial value $x_0$, satisfying
\begin{eqnarray*}
J(x_0, \bar{u}) \leq J(x_0, u), \quad \forall u\in L^2_{\mathcal{F}}(\mathbb{N}; \mathbb{R}^m),
\end{eqnarray*}
then, $\bar{u}$ is called as the open-loop optimal control with the initial value $x_0$;
Problem (SLQ) is (uniquely) open-loop solvable if it is (uniquely) open-loop solvable with any initial value $x_0$.
\end{definition}

\begin{definition}\label{D2}
Let $K: \mathbb{N} \rightarrow \mathbb{R}^{m\times n}$ be a  deterministic function and $v: \mathbb{N}\times \Omega \rightarrow \mathbb{R}^m$ be a $\mathcal{F}_{t-1}$-adapted process. ($K, v$) is called as a closed-loop strategy on $\mathbb{N}$ if $K\in L^2(\mathbb{N}; \mathbb{R}^{m\times n})$ and $v\in L^2_{\mathcal{F}}(\mathbb{N}; \mathbb{R}^{m})$, that is
$\sum\limits^{N-1}_{t=0} |K_t|^2< \infty$ and $E\Big(\sum\limits^{N-1}_{t=0} |v_t|^2\Big)< \infty$, and denote
$\mathcal{K}_{\mathbb{N}}$
%=\{(K_t, v_t)\Big|\sum\limits^{N-1}_{t=0} |K_t|^2< \infty, E\Big(\sum\limits^{N-1}_{t=0} |v_t|^2\Big)< \infty\}$
as the set of all closed-loop solutions;
$(K^*, v^*)\in \mathcal{K}_{\mathbb{N}}$ is called as the closed-loop optimal solution on $\mathbb{N}$ if
$\forall (K, v)\in\mathcal{K}_{\mathbb{N}}$, there holds
\begin{eqnarray*}
J(x_0, K^*x^*+v^*) \leq J(x_0, Kx+v),
\end{eqnarray*}
where $x^*$ and $x$ are the corresponding closed-loop systems satisfying
\begin{eqnarray*}
x_{t+1} &\hspace{-0.8em}=&\hspace{-0.8em} (A_t+B_tK_t)x_t+B_tv_t+b_t+[(C_t+D_tK_t)x_t\nonumber\\
&\hspace{-0.8em}&\hspace{-0.8em}+D_tv_t+\sigma_t]w_t,\\
x^*_{t+1} &\hspace{-0.8em}=&\hspace{-0.8em} (A_t+B_tK^*_t)x^*_t+B_tv^*_t+b_t
+[(C_t+D_tK^*_t)x^*_t\nonumber\\
&\hspace{-0.8em}&\hspace{-0.8em}+D_tv^*_t+\sigma_t]w_t,
\end{eqnarray*}
respectively, with the initial value $x_0$ and $x^*_0=x_0$;
If there exists a (uniquely) closed-loop optimal solution on $\mathbb{N}$, then Problem (SLQ) is called as (uniquely) closed-loop solvable.
\end{definition}

\begin{remark}
Following from the Definition \ref{D2}, $\forall x_0\in \mathbb{R}^n$, $(K^*, v^*)$ is an closed-loop optimal solution to Problem (SLQ) if and only if
\begin{eqnarray*}
J(x_0, K^*x^*+v^*) \leq J(x_0,u), \quad \forall u\in L^2_{\mathcal{F}}(\mathbb{N}; \mathbb{R}^m).
\end{eqnarray*}
Hence, it should be noted that $u^*=K^*x^*+v^*$ is also a open-loop optimal solution to Problem (SLQ) with the initial value $x^*_0$. That is to say, the existence of the closed-loop optimal solution to Problem (SLQ) indicates the existence of the open-loop optimal solution to it, but the reverse is not true, which can be seen in Example \ref{example1}.
When Problem (SLQ) is merely open-loop solvable, whether there also exists a control with a linear feedback form of the state is considered in this paper. To this end, the definition of a weak closed-loop solvability of Problem (SLQ) is given.
\end{remark}

\begin{definition}\label{definition3}
For any $\tilde{\mathbb{N}}\subseteq \mathbb{N}$ with $\tilde{\mathbb{N}}=\{0, 1, \cdots, m\}$ and $m\leq N-2$, let $K\in L^2(\tilde{\mathbb{N}}; \mathbb{R}^{m\times n})$ and
$v\in L^2_{\mathcal{F}}(\tilde{\mathbb{N}}; \mathbb{R}^{m})$. $(K_t, v_t)$ is called as  a weak closed-loop solution on $\tilde{\mathbb{N}}$ if for any initial value $x_0$, $u=Kx+v \in  L^2_{\mathcal{F}}(\mathbb{N}; \mathbb{R}^{m})$, where $x_t$ is the solution to the weak closed-loop system
\begin{eqnarray*}
x_{t+1} &\hspace{-0.8em}=&\hspace{-0.8em} (A_t+B_tK_t)x_t+B_tv_t+b_t+[(C_t+D_tK_t)x_t\nonumber\\
&\hspace{-0.8em}&\hspace{-0.8em}+D_tv_t+\sigma_t]w_t,
\end{eqnarray*}
with the initial value $x_0$, and denote $\mathcal{K}^w_{\mathbb{N}}$
%=\{(K_t, v_t)\Big|\sum\limits_{t\in \tilde{\mathbb{N}}} |K_t|^2< \infty, E\Big(\sum\limits_{t\in \tilde{\mathbb{N}}} |v_t|^2\Big)< \infty\}$
as the set of all weak closed-loop solutions; $(K^*_t, v^*_t)$ is said to be weak closed-loop optimal on $\tilde{\mathbb{N}}$, if for any initial value $x_0$ and $\forall (K, v)\in\mathcal{K}^w_{\mathbb{N}}$, there holds
\begin{eqnarray*}
J(x_0, K^*x^*+v^*) \leq J(x_0, Kx+v),
\end{eqnarray*}
where $x^*$ is the solution to the weak closed-loop system
\begin{eqnarray*}
x^*_{t+1} &\hspace{-0.8em}=&\hspace{-0.8em} (A_t+B_tK^*_t)x^*_t+B_tv^*_t+b_t+[(C_t+D_tK^*_t)x^*_t\nonumber\\
&\hspace{-0.8em}&\hspace{-0.8em}+D_tv^*_t+\sigma_t]w_t,
\end{eqnarray*}
with the initial value $x^*_0=x_0$;
Problem (SLQ) is said to be (uniquely) weakly closed-loop optimal solvable if there exists a (uniquely) weak closed-loop optimal solution $(K^*_t, v^*_t)$ on $\tilde{\mathbb{N}}$.
\end{definition}

\begin{remark}
Following from the Definition \ref{definition3}, $\forall x_0\in \mathbb{R}^n$, $(K^*, v^*)$ is a weakly closed-loop optimal solution to Problem (SLQ) if and only if
\begin{eqnarray*}
J(x_0, K^*x^*+v^*) \leq J(x_0,u), \quad \forall u\in L^2_{\mathcal{F}}({\mathbb{N}}; \mathbb{R}^{m}).
\end{eqnarray*}
\end{remark}

Based on the discussion above, the conditions of open-loop and closed-loop solvable for the Problem (SLQ) are derived.

\begin{lemma}\label{lemma1}
Problem (SLQ) is open-loop solvable if and only if $u_t$ subject to (\ref{1}) satisfies the equilibrium condition
\begin{eqnarray}\label{3}
  0 &\hspace{-0.8em}=&\hspace{-0.8em} \mathbb{E}[R_tu_t+S_tx_t+(B_t+w_tD_t)'\lambda_t+\rho_t|\mathcal{F}_{t-1}],
\end{eqnarray}
where the co-state with the terminal value $\lambda_{N-1}=Hx_{N}+g$  is such the following backward equation
\begin{eqnarray}\label{4}
  \lambda_{t-1}&\hspace{-0.8em}=&\hspace{-0.8em}\mathbb{E}[Q_tx_t+S'_tu_t+(A_t+w_tC_t)'\lambda_t+q_t|\mathcal{F}_{t-1}].
\end{eqnarray}
\end{lemma}

\emph{Proof:}
The proof can be obtained immediately from \cite{Zhang2015}.

\begin{lemma}\label{lemma2}
Let $(K^*, v^*)$ be a closed-loop control strategy of Problem (SLQ). Then $\forall x_0\in \mathbb{R}^n$, the following FBDEs admits an solution ($x^*_t, \lambda^*_{t-1}$) for $t\in \mathbb{N}$
\begin{eqnarray}
x^*_{t+1} &\hspace{-0.8em}=&\hspace{-0.8em} (A_t+B_tK^*_t)x^*_t+B_tv^*_t+b_t
+[(C_t+D_tK^*_t)x^*_t\nonumber\\
&\hspace{-0.8em}&\hspace{-0.8em}+D_tv^*_t+\sigma_t]w_t,\label{5-1}\\
\lambda^*_{t-1}&\hspace{-0.8em}=&\hspace{-0.8em}\mathbb{E}[(Q_t+S'_tK^*_t)x^*_t+S'_tv^*_t+(A_t+w_tC_t)\lambda^*_t
\nonumber\\
&\hspace{-0.8em}&\hspace{-0.8em}+q_t|\mathcal{F}_{t-1}],\label{5-2}
\end{eqnarray}
with the initial value $x^*_0=x_0$ and $\lambda^*_{N-1}=Hx^*_{N}+g$. Moreover, the equilibrium strategy satisfies
\begin{eqnarray}\label{6}
  0 &\hspace{-0.8em}=&\hspace{-0.8em} \mathbb{E}[(R_tK^*_t+S_t)x^*_t+R_tv^*_t+(B_t+w_tD_t)'\lambda^*_t+\rho_t|\mathcal{F}_{t-1}].
\end{eqnarray}
\end{lemma}

\emph{Proof:}
If $(K^*, v^*)$ is a closed-loop control strategy of Problem (SLQ), reconsidering (\ref{1})-(\ref{2}), one has
\begin{eqnarray*}
x_{t+1} &\hspace{-0.8em}=&\hspace{-0.8em} (A_t+B_tK^*_t)x_t+B_tv^*_t+b_t\nonumber\\
&\hspace{-0.8em}&\hspace{-0.8em}
+[(C_t+D_tK^*_t)x_t+D_tv^*_t+\sigma_t]w_t,\\
 J(x_0, K^*x^*+v^*) &\hspace{-0.8em}=&\hspace{-0.8em} \mathbb{E}\Big(\sum\limits^{N-1}_{t=0}\{x'_t[Q_t+(K^*_t)'R_tK^*_t+(K^*_t)'S_t\nonumber\\
&\hspace{-0.8em}&\hspace{-0.8em}+S'_tK^*_t]x_t
 +2(v^*_t)'[S_t+R_tK^*_t]x_t\nonumber\\
 &\hspace{-0.8em}&\hspace{-0.8em}+(v^*_t)'R_tv^*_t+2x'_tq_t+2(v^*_t)' \rho_t\nonumber\\
&\hspace{-0.8em}&\hspace{-0.8em}+2x'_t (K^*_t)'\rho_t\}+x'_{N}Hx_{N}+2x'_{N}g\Big).
\end{eqnarray*}
Combining with  Lemma \ref{lemma1},
(\ref{5-2}) and (\ref{6}) can be obtained accordingly.

%\begin{remark}
%$\forall x_0\in\mathbb{R}^n$, the FBDEs (\ref{5-1})-(\ref{6}) admits a solution and $(K^*_t, v^*_t)$ is independent of $x_0$.
%\end{remark}

The equivalent relationship between the closed-loop solvability of the Problem (SLQ) and the Riccati equation is given.
\begin{theorem}\label{theorem1}
Problem (SLQ) has a closed-loop solution $(K^*, v^*)$ if and only if the Riccati equation
with the terminal value $P_{N}=H$
\begin{eqnarray}\label{7}
P_t&\hspace{-0.8em}=&\hspace{-0.8em} Q_t+A'_tP_{t+1}A_t+C'_tP_{t+1}C_t-(A'_tP_{t+1}B_t+C'_tP_{t+1}D_t\nonumber\\
&\hspace{-0.8em}&\hspace{-0.8em}\times +S'_t)\hat{R}^{\dagger}_t(B'_tP_{t+1}A_t+D'_tP_{t+1}C_t+S_t),
\end{eqnarray}
admits a solution $P_t$ such that
\begin{eqnarray}
&\hspace{-0.8em}&\hspace{-0.8em}\hat{R}_t \geq 0,\label{8-1}\\
&\hspace{-0.8em}&\hspace{-0.8em}\hat{K}_t=- \hat{R}^{\dagger}_t(B'_tP_{t+1}A_t+D'_tP_{t+1}C_t+S_t)
\in L^2(\mathbb{N}; \mathbb{R}^{m\times n}),\label{8-2}\\
&\hspace{-0.8em}&\hspace{-0.8em}Range(B'_tP_{t+1}A_t+D'_tP_{t+1}C_t+S_t)  \subseteq  Range( \hat{R}_t),\label{8-3}
\end{eqnarray}
where $\hat{R}_t={R}_t+B'_tP_{t+1}B_t+D'_tP_{t+1}D_t$, and the adapted solution ($\eta_{t-1}, v_t$) of BSDE
\begin{eqnarray}\label{8-4}
 \eta_{t-1}&\hspace{-0.8em}=&\hspace{-0.8em}\mathbb{E}\{[(A_t+w_tC_t)'+\hat{K}'_t(B_t+w_tD_t)']
 \eta_t|\mathcal{F}_{t-1}\}\nonumber\\
 &\hspace{-0.8em}&\hspace{-0.8em}+(A'_t+\hat{K}'_tB'_t)P_{t+1}b_t+(C'_t+\hat{K}'_tD'_t)P_{t+1}\sigma_t
 +q_t\rho_t\nonumber\\
 &\hspace{-0.8em}&\hspace{-0.8em}+\hat{K}'_t
\end{eqnarray}
with the terminal value $\eta_{N-1}=g$ satisfies
\begin{eqnarray}
&\hspace{-0.8em}&\hspace{-0.8em}\mathbb{E}[(B_t+w_tD_t)'\eta_t|\mathcal{F}_{t-1}]+B'_tP_{t+1}b_t
+D'_tP_{t+1}\sigma_t\nonumber\\
&\hspace{-0.8em}&\hspace{-0.8em}+\rho_t\in Range( \hat{R}_t),\label{8-5}\\
\hat{v}_t&\hspace{-0.8em}=&\hspace{-0.8em}-\hat{R}_t^{\dagger}\{\mathbb{E}[(B_t+w_tD_t)'\eta_t|
\mathcal{F}_{t-1}]+B'_tP_{t+1}b_t+D'_tP_{t+1}\sigma_t\nonumber\\
&\hspace{-0.8em}&\hspace{-0.8em}+\rho_t\}\in L^2_{\mathcal{F}}(\mathbb{N}; \mathbb{R}^{m})\label{8-6}
\end{eqnarray}

In this case, the closed-loop solution $(K^*, v^*)$ is given by
\begin{eqnarray}
K^*_t&\hspace{-0.8em}=&\hspace{-0.8em} \hat{K}_t+(I-\hat{R}_t^{\dagger}\hat{R}_t) z_t,\label{9-1}\\
v^*_t&\hspace{-0.8em}=&\hspace{-0.8em}\hat{v}_t+(I-\hat{R}_t^{\dagger}\hat{R}_t) y_t\label{9-2},
%&\hspace{-0.8em}=&\hspace{-0.8em}[I-(R_t+B'_tP_{t+1}B_t)^{\dagger}(R_t+B'_tP_{t+1}B_t)]y_t,{\color{blue}?????}\label{9-3}
\end{eqnarray}
where $z\in L^2(\mathbb{N}; \mathbb{R}^{m\times n})$ and $y\in L^2_{\mathcal{F}}(\mathbb{N}; \mathbb{R}^{m})$ are arbitrary.

Moreover, the value function is such that
\begin{eqnarray}\label{10}
V(x_0)&\hspace{-0.8em}=&\hspace{-0.8em}\mathbb{E}(x'_0P_0x_0+2x'_0\eta_{-1})
+\sum^{N-1}_{t=0}E\{b'_tP_{t+1}b_t+\sigma'_tP_{t+1} \nonumber\\
&\hspace{-0.8em}&\hspace{-0.8em}\times \sigma_t+2(b_t+w_t\sigma_t)'\eta_t
-[\mathbb{E}((B_t+w_tD_t)'\eta_{t}|\mathcal{F}_{t-1})\nonumber\\
&\hspace{-0.8em}&\hspace{-0.8em}+B'_tP_{t+1}b_t+D'_tP_{t+1}\sigma_t
+\rho_t]'
\hat{R}_t[\mathbb{E}((B_t+w_tD_t)'\nonumber\\
&\hspace{-0.8em}&\hspace{-0.8em}\times \eta_{t}|\mathcal{F}_{t-1})+B'_tP_{t+1}b_t+D'_tP_{t+1}\sigma_t+\rho_t]\}.
\end{eqnarray}
\end{theorem}

\emph{Proof:}
``Necessary.'' Suppose $(K^*, v^*)$ is the optimal closed-loop solution of the Problem (SLQ), according to Lemma \ref{lemma2}, $(K^*, v^*)$ satisfies FBDEs (\ref{5-1})-(\ref{6}) for any $x_0\in \mathbb{R}^n$.
Let ($\bar{x}_t, \bar{\lambda}_{t-1}$) on $\mathbb{N}$ be the solution to the FBDEs
\begin{eqnarray}
\bar{x}_{t+1} &\hspace{-0.8em}=&\hspace{-0.8em} (A_t+B_tK^*_t)\bar{x}_t+(C_t+D_tK^*_t)w_t\bar{x}_t, \bar{x}_0=x_0\label{11-1}\\
\bar{\lambda}_{t-1}&\hspace{-0.8em}=&\hspace{-0.8em}(Q_t+S'_tK^*_t)\bar{x}_t+\mathbb{E}[(A_t+w_tC_t)'
\bar{\lambda}_t|\mathcal{F}_{t-1}],\label{11-2}\\ \bar{\lambda}_{N-1}&\hspace{-0.8em}=&\hspace{-0.8em}H\bar{x}_{N},\nonumber
\end{eqnarray}
then, there exists the following equilibrium condition
\begin{eqnarray}\label{12}
  0 &\hspace{-0.8em}=&\hspace{-0.8em} (R_tK^*_t+S_t)\bar{x}_t+\mathbb{E}[(B_t+w_tD_t)'\bar{\lambda}_t|\mathcal{F}_{t-1}].
\end{eqnarray}

Solving (\ref{11-1})-(\ref{11-2}) by backward iteration from $t=N-1$, one has
\begin{eqnarray*}
\bar{\lambda}_{t-1}&\hspace{-0.8em}=&\hspace{-0.8em}(Q_t+S'_tK^*_t)\bar{x}_t+\mathbb{E}
\{(A_t+w_tC_t)'P_{t+1}\nonumber\\
&\hspace{-0.8em}&\hspace{-0.8em}\times[(A_t+B_tK^*_t)+(C_t+D_tK^*_t)w_t]\bar{x}_t|\mathcal{F}_{t-1}\}\nonumber\\
&\hspace{-0.8em}=&\hspace{-0.8em}[Q_t+A'_tP_{t+1}A_t+C'_tP_{t+1}C_t+(A'_tP_{t+1}B_t\nonumber\\
&\hspace{-0.8em}&\hspace{-0.8em}+C'_tP_{t+1}D_t+S'_t)K^*_t]\bar{x}_t\nonumber\\
&\doteq& P_t\bar{x}_t,
\end{eqnarray*}
with the terminal value $P_{N}=H$, i.e.,
\begin{eqnarray}\label{12-1}
  P_t &\hspace{-0.8em}=&\hspace{-0.8em} Q_t+A'_tP_{t+1}A_t+C'_tP_{t+1}C_t\nonumber\\
&\hspace{-0.8em}&\hspace{-0.8em}+(A'_tP_{t+1}B_t+C'_tP_{t+1}D_t+S'_t)K^*_t.
\end{eqnarray}
Adding $\bar{\lambda}_{t}=P_{t+1}\bar{x}_{t+1}$ into (\ref{12}), it yields that
\begin{eqnarray*}
0 &\hspace{-0.8em}=&\hspace{-0.8em} (R_tK^*_t+S_t)\bar{x}_t+\mathbb{E}\{(B_t+w_tD_t)'P_{t+1}[(A_t+B_tK^*_t)\nonumber\\
&\hspace{-0.8em}&\hspace{-0.8em}
+(C_t+D_tK^*_t)w_t]\bar{x}_t|\mathcal{F}_{t-1}\}\nonumber\\
  &\hspace{-0.8em}=&\hspace{-0.8em}(\hat{R}_tK^*_t+B'_tP_{t+1}A_t+D'_tP_{t+1}C_t+S_t)\bar{x}_t.
\end{eqnarray*}
The arbitrariness of $\bar{x}_t$ means that
\begin{eqnarray}\label{13}
\hat{R}_tK^*_t+B'_tP_{t+1}A_t+D'_tP_{t+1}C_t+S_t=0,
\end{eqnarray}
which implies $Range(B'_tP_{t+1}A_t+D'_tP_{t+1}C_t+S_t)  \subseteq  Range( \hat{R}_t)$.
To this end, $P_t$ can be rewritten as
\begin{eqnarray}\label{12-2}
P_t&\hspace{-0.8em}=&\hspace{-0.8em}Q_t+A'_tP_{t+1}A_t+C'_tP_{t+1}C_t+(A'_tP_{t+1}B_t\nonumber\\
&\hspace{-0.8em}&\hspace{-0.8em}+C'_tP_{t+1}D_t+S'_t)K^*_t+(K^*_t)'[\hat{R}_tK^*_t+B'_tP_{t+1}A_t\nonumber\\
&\hspace{-0.8em}&\hspace{-0.8em}+D'_tP_{t+1}C_t+S_t]\nonumber\\
&\hspace{-0.8em}=&\hspace{-0.8em}Q_t+(A_t+B_tK^*_t)'P_{t+1}(A_t+B_tK^*_t)+(C_t+D_tK^*_t)'\nonumber\\
&\hspace{-0.8em}&\hspace{-0.8em}\times P_{t+1}(C_t+D_tK^*_t)
+S'_tK^*_t+(K^*_t)'S_t+(K^*_t)'R_tK^*_t.\nonumber\\
\end{eqnarray}

Following from (\ref{13}), one has
\begin{eqnarray*}
 \hat{R}_t^{\dagger}(B'_tP_{t+1}A_t+D'_tP_{t+1}C_t+S_t)=-\hat{R}_t^{\dagger} \hat{R}_tK^*_t,
\end{eqnarray*}
that is to say
\begin{eqnarray*}
K^*_t&\hspace{-0.8em}=&\hspace{-0.8em} - R_t^{\dagger}(B'_tP_{t+1}A_t+D'_tP_{t+1}B_t+S_t)+(I-\hat{R}_t^{\dagger}\hat{R}_t) z_t.
\end{eqnarray*}

Consequently,
\begin{eqnarray*}
&\hspace{-0.8em}&\hspace{-0.8em}(A'_tP_{t+1}B_t+C'_tP_{t+1}D_t+S'_t)K^*_t\nonumber\\
=&\hspace{-0.8em}&\hspace{-0.8em} -(K^*_t)'\hat{R}_tK^*_t \nonumber\\
=&\hspace{-0.8em}&\hspace{-0.8em}(K^*_t)'\hat{R}_t\hat{R}_t^{\dagger}(B'_tP_{t+1}A_t+D'_tP_{t+1}C_t+S_t)
\nonumber\\
=&\hspace{-0.8em}&\hspace{-0.8em}-(A'_tP_{t+1}B_t+C'_tP_{t+1}D_t+S'_t)\hat{R}_t^{\dagger}
(B'_tP_{t+1}A_t+D'_tP_{t+1}C_t\nonumber\\
&\hspace{-0.8em}&\hspace{-0.8em}+S_t).
\end{eqnarray*}
Substituting the above equation to (\ref{12-1}), the Riccati equation (\ref{7}) is obtained.

Next, we will calculate $v^*_t$ on $\mathbb{N}$. Denote $\eta_{t-1}=\lambda^*_{t-1}-P_tx^*_t$ with the terminal value $\eta_{N-1}=g$. Adding it into (\ref{6}), it derives
\begin{eqnarray}\label{13-1}
  0 &\hspace{-0.8em}=&\hspace{-0.8em} (R_tK^*_t+S_t)x^*_t+R_tv^*_t+\mathbb{E}[(B_t+w_tD_t)'(\eta_{t}+P_{t+1}\nonumber\\
&\hspace{-0.8em}&\hspace{-0.8em}\times x^*_{t+1})
+\rho_t|\mathcal{F}_{t-1}] (R_tK^*_t+S_t)x^*_t+R_tv^*_t+\mathbb{E}[(B_t\nonumber\\
&\hspace{-0.8em}&\hspace{-0.8em}+w_tD_t)'\eta_{t}|\mathcal{F}_{t-1}]
+\mathbb{E}\{(B_t +w_tD_t)'P_{t+1}[(A_t+B_tK^*_t)\nonumber\\
&\hspace{-0.8em}&\hspace{-0.8em}\times x^*_t+B_tv^*_t+b_t+((C_t+D_tK^*_t)x^*_t
+D_tv^*_t+\sigma_t)\nonumber\\
&\hspace{-0.8em}&\hspace{-0.8em}\times w_t]|\mathcal{F}_{t-1}\}
+\rho_t\nonumber\\
  &\hspace{-0.8em}=&\hspace{-0.8em}\hat{R}_tv^*_t
+[B'_tP_{t+1}A_t+D'_tP_{t+1}C_t+S_t+\hat{R}_tK^*_t]x^*_t\nonumber\\
&\hspace{-0.8em}&\hspace{-0.8em}
+\mathbb{E}[(B_t+w_tD_t)'\eta_{t}|\mathcal{F}_{t-1}]+B'_tP_{t+1}b_t+D'_tP_{t+1}\sigma_t+\rho_t\nonumber\\
&\hspace{-0.8em}=&\hspace{-0.8em}\hat{R}_tv^*_t+\mathbb{E}[(B_t+w_tD_t)'\eta_{t}|\mathcal{F}_{t-1}]
+B'_tP_{t+1}b_t\nonumber\\
&\hspace{-0.8em}&\hspace{-0.8em}+D'_tP_{t+1}\sigma_t+\rho_t,
\end{eqnarray}
which implies (\ref{8-5}) is established.

On the other hand, following from (\ref{13-1}), one has
\begin{eqnarray*}
\hat{R}_t^{\dag}\hat{R}_tv^*_t&\hspace{-0.8em}=&\hspace{-0.8em}-\hat{R}_t^{\dag} \{\mathbb{E}[(B_t+w_tD_t)'\eta_{t}|\mathcal{F}_{t-1}]
+B'_tP_{t+1}b_t\nonumber\\
&\hspace{-0.8em}&\hspace{-0.8em}+D'_tP_{t+1}\sigma_t+\rho_t\}.
\end{eqnarray*}
Since $\hat{R}_t^{\dag}\hat{R}_t$ is an orthogonal projection, (\ref{8-6}) is established and $v^*_t$ is written as (\ref{9-2}) for $t\in\mathbb{N}$.
%\begin{eqnarray}\label{14}
%% \nonumber to remove numbering (before each equation)
%v^*_t&\hspace{-0.8em}=&\hspace{-0.8em}-(R_t+B'_tP_{t+1}B_t)^{\dagger}B'_t\eta_{t}
%+[I-(R_t+B'_tP_{t+1}B_t)^{\dagger}(R_t+B'_tP_{t+1}B_t)]y_t,
%\end{eqnarray}
%where $y_t\in L^2(0, N; \mathbb{R}^m)$ is arbitrary. And

According to (\ref{9-2}), we have
\begin{eqnarray*}
&\hspace{-0.8em}&\hspace{-0.8em}(A'_tP_{t+1}B_t+C'_tP_{t+1}D_t+S'_t)v^*_t\nonumber\\
=&\hspace{-0.8em}&\hspace{-0.8em}-(A'_tP_{t+1}B_t+C'_tP_{t+1}D_t+S'_t)\hat{R}_t^{\dagger} \{\mathbb{E}[(B_t+w_tD_t)'\nonumber\\
&\hspace{-0.8em}&\hspace{-0.8em}\times \eta_{t}|\mathcal{F}_{t-1}]+B'_tP_{t+1}b_t+D'_tP_{t+1}\sigma_t
+\rho_t\}
+(A'_tP_{t+1}B_t\nonumber\\
&\hspace{-0.8em}&\hspace{-0.8em}+C'_tP_{t+1}D_t+S'_t)(I-\hat{R}_t^{\dagger} \hat{R}_t)y_t\nonumber\\
=&\hspace{-0.8em}&\hspace{-0.8em}-(A'_tP_{t+1}B_t+C'_tP_{t+1}D_t+S'_t)\hat{R}_t^{\dagger}
\{\mathbb{E}[(B_t+w_tD_t)'\nonumber\\
&\hspace{-0.8em}&\hspace{-0.8em}\times \eta_{t}|\mathcal{F}_{t-1}]+B'_tP_{t+1}b_t+D'_tP_{t+1}\sigma_t+\rho_t\}
  -(K^*_t)'\hat{R}_t\nonumber\\
&\hspace{-0.8em}&\hspace{-0.8em}\times (I-\hat{R}_t^{\dagger}\hat{R}_t)y_t\nonumber\\
=&\hspace{-0.8em}&\hspace{-0.8em}-(A'_tP_{t+1}B_t+C'_tP_{t+1}D_t+S'_t)\hat{R}_t^{\dagger}
\{\mathbb{E}[(B_t+w_tD_t)'\nonumber\\
&\hspace{-0.8em}&\hspace{-0.8em}\times \eta_{t}|\mathcal{F}_{t-1}]+B'_tP_{t+1}b_t+D'_tP_{t+1}\sigma_t+\rho_t\}\nonumber\\
=&\hspace{-0.8em}&\hspace{-0.8em}-\hat{K}'_t\{\mathbb{E}[(B_t+w_tD_t)'\eta_{t}|\mathcal{F}_{t-1}]
+B'_tP_{t+1}b_t+D'_tP_{t+1}\sigma_t\nonumber\\
&\hspace{-0.8em}&\hspace{-0.8em}+\rho_t\}.
\end{eqnarray*}

Combining with (\ref{5-2}), $\eta_{t-1}$ is calculated as
\begin{eqnarray}\label{15}
  \eta_{t-1} &\hspace{-0.8em}=&\hspace{-0.8em} (Q_t+S'_tK^*_t)x^*_t+S'_tv^*_t+\mathbb{E}[(A_t+w_tC_t)'(\eta_t+P_{t+1}\nonumber\\
&\hspace{-0.8em}&\hspace{-0.8em}\times x^*_{t+1})|\mathcal{F}_{t-1}]
  +q_t-P_tx^*_t \nonumber\\
&\hspace{-0.8em}=&\hspace{-0.8em}(Q_t+S'_tK^*_t)x^*_t+S'_tv^*_t+\mathbb{E}[(A_t+w_tC_t)'\eta_t
|\mathcal{F}_{t-1}]\nonumber\\
&\hspace{-0.8em}&\hspace{-0.8em}
+\mathbb{E}\{(A_t+w_tC_t)'P_{t+1} [(A_t+B_tK^*_t)x^*_t+B_tv^*_t+b_t\nonumber\\
&\hspace{-0.8em}&\hspace{-0.8em}+((C_t+D_tK^*_t)x^*_t
+D_tv^*_t+\sigma_t)w_t]|\mathcal{F}_{t-1}\}+q_t \nonumber\\
&\hspace{-0.8em}=&\hspace{-0.8em}-P_tx^*_t(A'_tP_{t+1}B_t+C'_tP_{t+1}D_t+S'_t)v^*_t\nonumber\\
&\hspace{-0.8em}&\hspace{-0.8em}+\mathbb{E}
[(A_t+w_tC_t)'\eta_t|\mathcal{F}_{t-1}]
+[Q_t+A'_tP_{t+1}A_t\nonumber\\
&\hspace{-0.8em}&\hspace{-0.8em}+C'_tP_{t+1}C_t+(A'_tP_{t+1}B_t+C'_tP_{t+1}D_t+S'_t)K^*_t\nonumber\\
&\hspace{-0.8em}&\hspace{-0.8em}-P_t]x^*_t
+A'_tP_{t+1}b_t+C'_tP_{t+1}\sigma_t+q_t.\nonumber\\
\end{eqnarray}

According to (\ref{9-2}), we have
\begin{eqnarray*}
&\hspace{-0.8em}&\hspace{-0.8em}(A'_tP_{t+1}B_t+C'_tP_{t+1}D_t+S'_t)v^*_t\nonumber\\
=&\hspace{-0.8em}&\hspace{-0.8em}-(A'_tP_{t+1}B_t+C'_tP_{t+1}D_t+S'_t)\hat{R}_t^{\dagger} \{\mathbb{E}[(B_t+w_tD_t)'\nonumber\\
&\hspace{-0.8em}&\hspace{-0.8em}\times \eta_{t}|\mathcal{F}_{t-1}]+B'_tP_{t+1}b_t+D'_tP_{t+1}\sigma_t
+\rho_t\}\nonumber\\
&\hspace{-0.8em}&\hspace{-0.8em}
+(A'_tP_{t+1}B_t+C'_tP_{t+1}D_t+S'_t)(I-\hat{R}_t^{\dagger} \hat{R}_t)y_t\nonumber\\
=&\hspace{-0.8em}&\hspace{-0.8em}-(A'_tP_{t+1}B_t+C'_tP_{t+1}D_t+S'_t)\hat{R}_t^{\dagger}
\{\mathbb{E}[(B_t+w_tD_t)'\nonumber\\
&\hspace{-0.8em}&\hspace{-0.8em}\times \eta_{t}|\mathcal{F}_{t-1}]+B'_tP_{t+1}b_t+D'_tP_{t+1}\sigma_t+\rho_t\}
  -(K^*_t)'\hat{R}_t\nonumber\\
&\hspace{-0.8em}&\hspace{-0.8em}\times (I-\hat{R}_t^{\dagger}\hat{R}_t)y_t\nonumber\\
=&\hspace{-0.8em}&\hspace{-0.8em}-(A'_tP_{t+1}B_t+C'_tP_{t+1}D_t+S'_t)\hat{R}_t^{\dagger}
\{\mathbb{E}[(B_t+w_tD_t)'\nonumber\\
&\hspace{-0.8em}&\hspace{-0.8em}\times \eta_{t}|\mathcal{F}_{t-1}]+B'_tP_{t+1}b_t+D'_tP_{t+1}\sigma_t+\rho_t\}\nonumber\\
=&\hspace{-0.8em}&\hspace{-0.8em}-\hat{K}'_t\{\mathbb{E}[(B_t+w_tD_t)'\eta_{t}|\mathcal{F}_{t-1}]\nonumber\\
&\hspace{-0.8em}&\hspace{-0.8em}+B'_tP_{t+1}b_t+D'_tP_{t+1}\sigma_t+\rho_t\}.
\end{eqnarray*}

Adding the above equation into (\ref{15}), $\eta_{t-1}$ degrades  into
\begin{eqnarray*}
\eta_{t-1} &\hspace{-0.8em}=&\hspace{-0.8em} \hat{K}'_t\{\mathbb{E}[(B_t+w_tD_t)'\eta_{t}|\mathcal{F}_{t-1}]
+B'_tP_{t+1}b_t+D'_tP_{t+1}\nonumber\\
&\hspace{-0.8em}&\hspace{-0.8em}\times \sigma_t+\rho_t\}+\mathbb{E}[(A_t+w_tC_t)'\eta_t|\mathcal{F}_{t-1}]
+[Q_t\nonumber\\
&\hspace{-0.8em}&\hspace{-0.8em}+A'_tP_{t+1}A_t+C'_tP_{t+1}C_t+(A'_tP_{t+1}B_t
+C'_tP_{t+1}D_t\nonumber\\
&\hspace{-0.8em}&\hspace{-0.8em}+S'_t)K^*_t-P_t]x^*_t
+A'_tP_{t+1}b_t+C'_tP_{t+1}\sigma_t+q_t\nonumber\\
&\hspace{-0.8em}=&\hspace{-0.8em}\mathbb{E}\{[(A_t+w_tC_t)'+\hat{K}'_t(B_t+w_tD_t)']\eta_{t}
|\mathcal{F}_{t-1}\}\nonumber\\
&\hspace{-0.8em}&\hspace{-0.8em}
+(A_t+B_t\hat{K}_t)'P_{t+1}b_t+(C_t+D_t\hat{K}_t)'P_{t+1}\sigma_t\nonumber\\
&\hspace{-0.8em}&\hspace{-0.8em}+\hat{K}'_t\rho_t+q_t,
\end{eqnarray*}
where the last equality is established due to (\ref{12-1}).

``Sufficiency''. $ \forall u\in L^2_{\mathcal{F}}({\mathbb{N}}; \mathbb{R}^{m})$, let $x_t=x(t; \mathbb{N}, u_t)$ is the corresponding state process. From (\ref{1}) and (\ref{7}), it derives
\begin{eqnarray}\label{15-1}
&\hspace{-0.8em}&\hspace{-0.8em}\mathbb{E}[x'_tP_tx_t-x'_{t+1}P_{t+1}x_{t+1}]\nonumber\\
=&\hspace{-0.8em}&\hspace{-0.8em}\mathbb{E}\{x'_t(P_t-Q_t-A'_tP_{t+1}A_t-C'_tP_{t+1}C_t)x_t
-u'_t\hat{R}_tu_t\nonumber\\
&\hspace{-0.8em}&\hspace{-0.8em}-2u'_t[(B'_tP_{t+1}A_t+D'_tP_{t+1}C_t+S_t)x_t+B'_tP_{t+1}b_t\nonumber\\
&\hspace{-0.8em}&\hspace{-0.8em}+D'_tP_{t+1}\sigma_t+\rho_t]
-2x'_t[A'_tP_{t+1}b_t+C'_tP_{t+1}\sigma_t+q_t]\nonumber\\
&\hspace{-0.8em}&\hspace{-0.8em}-b'_tP_{t+1}b_t-\sigma'_tP_{t+1}\sigma_t
+(x'_tQ_tx_t+2u'_tS_tx_t+u'_tR_tu_t\nonumber\\
&\hspace{-0.8em}&\hspace{-0.8em}+2x'_tq_t+2u'_t \rho_t)\}.
\end{eqnarray}
And according to (\ref{1}) and (\ref{8-4}), one has
\begin{eqnarray}\label{15-2}
\mathbb{E}[x'_t\eta_{t-1}-x'_{t+1}\eta_{t}]&\hspace{-0.8em}=&\hspace{-0.8em}-2\mathbb{E}\{x'_t
(B'_tP_{t+1}A_t+D'_tP_{t+1}C_t\nonumber\\
&\hspace{-0.8em}&\hspace{-0.8em}+S_t)'
\hat{R}_t^{\dag}[(B_t+w_tD_t)'\eta_{t}+B'_tP_{t+1}b_t\nonumber\\
&\hspace{-0.8em}&\hspace{-0.8em}+D'_tP_{t+1}\sigma_t+\rho_t]-u'_t(B_t+w_tD_t)'\eta_{t}\nonumber\\
&\hspace{-0.8em}&\hspace{-0.8em}
-(b_t+w_t\sigma_t)'\eta_{t}\}.
\end{eqnarray}
Combining (\ref{15-1}) with (\ref{15-2}), there holds that
\begin{eqnarray}\label{15-3-1}
&\hspace{-0.8em}&\hspace{-0.8em}\mathbb{E}[x'_tP_tx_t-x'_{t+1}P_{t+1}x_{t+1}+x'_t\eta_{t-1}-x'_{t+1}\eta_{t}]\nonumber\\
=&\hspace{-0.8em}&\hspace{-0.8em}\mathbb{E}\{-x'_t(A'_tP_{t+1}B_t+C'_tP_{t+1}D_t+S'_t)
\hat{R}_t^{\dagger}(B'_tP_{t+1}A_t\nonumber\\
&\hspace{-0.8em}&\hspace{-0.8em}+D'_tP_{t+1}C_t+S_t)x_t
-u'_t\hat{R}_tu_t-2u'_t[(B'_tP_{t+1}A_t\nonumber\\
&\hspace{-0.8em}&\hspace{-0.8em}+D'_tP_{t+1}C_t+S_t)x_t+(B_t+w_tD_t)'\eta_{t}+B'_tP_{t+1}b_t\nonumber\\
&\hspace{-0.8em}&\hspace{-0.8em}
+D'_tP_{t+1}\sigma_t+\rho_t]-2x'_t(B'_tP_{t+1}A_t+D'_tP_{t+1}C_t+S_t)'\nonumber\\
&\hspace{-0.8em}&\hspace{-0.8em}\times
\hat{R}_t^{\dag}[(B_t+w_tD_t)'\eta_{t}
+B'_tP_{t+1}b_t+D'_tP_{t+1}\sigma_t+\rho_t]\nonumber\\
&\hspace{-0.8em}&\hspace{-0.8em}
-b'_tP_{t+1}b_t-\sigma'_tP_{t+1}\sigma_t-2(b_t+w_t\sigma_t)'\eta_{t}
+(x'_tQ_tx_t\nonumber\\
&\hspace{-0.8em}&\hspace{-0.8em}+2u'_tS_tx_t+u'_tR_tu_t+2x'_tq_t+2u'_t \rho_t)\}.
\end{eqnarray}
By using (\ref{13}) and (\ref{13-1}), (\ref{15-3-1}) is rewritten as
\begin{eqnarray}\label{15-3}
&\hspace{-0.8em}&\hspace{-0.8em}\mathbb{E}[x'_tP_tx_t-x'_{t+1}P_{t+1}x_{t+1}+x'_t\eta_{t-1}-x'_{t+1}\eta_{t}]\nonumber\\
=&\hspace{-0.8em}&\hspace{-0.8em}\mathbb{E}\{-x'_t(K^*_t)'\hat{R}_t
\hat{R}_t^{\dagger}\hat{R}_tK^*_tx_t
-u'_t\hat{R}_tu_t+2u'_t\hat{R}_t(K^*_tx_t\nonumber\\
&\hspace{-0.8em}&\hspace{-0.8em}+v^*_t)
-2x'_t(K^*_t)'\hat{R}_t
\hat{R}_t^{\dagger}\hat{R}_tv^*_t
-b'_tP_{t+1}b_t-\sigma'_tP_{t+1}\sigma_t\nonumber\\
&\hspace{-0.8em}&\hspace{-0.8em}-2(b_t+w_t\sigma_t)'\eta_{t}
+(x'_tQ_tx_t+2u'_tS_tx_t+u'_tR_tu_t\nonumber\\
&\hspace{-0.8em}&\hspace{-0.8em}+2x'_tq_t+2u'_t \rho_t)\}\nonumber\\
=&\hspace{-0.8em}&\hspace{-0.8em}\mathbb{E}\{-[u_t-K^*_tx_t-v^*_t]'\hat{R}_t[u_t-K^*_tx_t-v^*_t]
+[\mathbb{E}((B_t\nonumber\\
&\hspace{-0.8em}&\hspace{-0.8em}+w_tD_t)'\eta_{t}|\mathcal{F}_{t-1})+B'_tP_{t+1}b_t
+D'_tP_{t+1}\sigma_t+\rho_t]'
\hat{R}_t\nonumber\\
&\hspace{-0.8em}&\hspace{-0.8em}\times [\mathbb{E}((B_t+w_tD_t)'\eta_{t}|\mathcal{F}_{t-1})
+B'_tP_{t+1}b_t+D'_tP_{t+1}\sigma_t\nonumber\\
&\hspace{-0.8em}&\hspace{-0.8em}+\rho_t]-b'_tP_{t+1}b_t-\sigma'_tP_{t+1}\sigma_t-2(b_t+w_t\sigma_t)'
\eta_{t}\nonumber\\
&\hspace{-0.8em}&\hspace{-0.8em}
+(x'_tQ_tx_t+2u'_tS_tx_t+u'_tR_tu_t+2x'_tq_t+2u'_t \rho_t)\}.
\end{eqnarray}

Substituting from $k=0$ to $k=N-1$ on both sides of (\ref{15-3}), and following from (\ref{2}), the cost function is re-calculated as
\begin{eqnarray}\label{16}
J(x_0, u) &\hspace{-0.8em}=&\hspace{-0.8em} \mathbb{E}[(u_t-K^*_tx_t-v^*_t)'\hat{R}_t(u_t-K^*_tx_t-v^*_t)]\nonumber\\
&\hspace{-0.8em}&\hspace{-0.8em}+J(x_0, K^*x^*+v^*).
\end{eqnarray}
Based on the discussion above, it can be concluded that
\begin{eqnarray*}
J(x_0, K^*x^*+v^*) \leq J(x_0, Kx+v)
\end{eqnarray*}
if and only if $\hat{R}_t\geq 0$ for $t\in \mathbb{N}$.  That is, $(K^*, v^*)$ is the closed-loop optimal solution to  Problem (SLQ). This completes the proof.

\begin{remark}
The closed-loop solution $(K^*, v^*)$ is independent of the initial value $x_0$.
In Theorem \ref{theorem1}, it has shown that  Problem (SLQ) is closed-loop optimal solvable, which implies the open-loop solvability of Problem (SLQ),  if and only if the Riccati equation (\ref{7}) admits a regular solution.
Thus, in Example \ref{example1}, the closed-loop solution calculated from Riccati equation (\ref{7}) cannot be find due to
\begin{eqnarray}
 Range( R+B'P_{t+1}B+D'P_{t+1}D)&\hspace{-0.8em}=&\hspace{-0.8em}Range( 0)=\{0\},\\
 Range(B'P_{t+1}A+D'P_{t+1}C+S)&\hspace{-0.8em}=&\hspace{-0.8em}Range(1)=\mathbb{R}.
\end{eqnarray}
That is,
\begin{eqnarray*}
Range(B'P_{t+1}A+D'P_{t+1}C+S] &\hspace{-0.8em} \nsubseteq&\hspace{-0.8em}  Range( R+B'P_{t+1}B\nonumber\\
&\hspace{-0.8em}&\hspace{-0.8em}+D'P_{t+1}D].
\end{eqnarray*}
\end{remark}

The relationship between the convexity of the map $u_t\mapsto J^0(x_0, u)$ and the regular solvability of the Riccati equation (\ref{7}) is given.
\begin{theorem}\label{theorem2}
(1) Suppose  Problem (SLQ) is open-loop solvable, then, for any $u\in L^2_{\mathcal{F}}({\mathbb{N}}; \mathbb{R}^{m})$, $J^0(0, u)\geq 0$.

(2) If $J^0(0, u)\geq 0$, then,  the map $u_t\mapsto J(x_0, u)$ is convex.

(3) Suppose there exists a constant $\alpha >0$ such that
\begin{eqnarray*}
 J^0(0, u)\geq \alpha \mathbb{E} \Big(\sum\limits^{N-1}_{t=0}|u_t|^2\Big), \quad \forall u\in L^2_{\mathcal{F}}({\mathbb{N}}; \mathbb{R}^{m}),
\end{eqnarray*}
Then, the Riccati equation (\ref{7}) has a unique solution $P_t$ satisfying
\begin{eqnarray*}
\hat{R}_t \geq \alpha I.
\end{eqnarray*}
To this end, Problem (SLQ) is uniquely closed-loop solvable, and it's uniquely open-loop solvable. The unique closed-loop optimal solution is given by
\begin{eqnarray*}
K^*_t&\hspace{-0.8em}=&\hspace{-0.8em} - \hat{R}_t^{-1}(B'_tP_{t+1}A_t++D'_tP_{t+1}C_t+S_t),\\
v^*_t&\hspace{-0.8em}=&\hspace{-0.8em}-\hat{R}_t^{-1}\{\mathbb{E}[(B_t+w_tD_t)'\eta_t|
\mathcal{F}_{t-1}]+B'_tP_{t+1}b_t\nonumber\\
&\hspace{-0.8em}&\hspace{-0.8em}+D'_tP_{t+1}\sigma_t+\rho_t\},
\end{eqnarray*}
where $\eta_t$ is the adapted solution to BSDE (\ref{8-4}).
Accordingly, the open-loop optimal control of the Problem (SLQ) with initial value $x_0$ is such that $u^*_t=K^*_tx^*_t+v^*_t$ on $\mathbb{N}$, where $x^*_t$ satisfies the closed-loop system (\ref{5-1}) with initial value $x^*_0=x_0$.
\end{theorem}

\emph{Proof:}
(1) Suppose $\tilde{x}_t$ is the solution for (\ref{2-1}) on $\mathbb{N}$ with the initial value ${x}_0=0$ for all $ u\in L^2_{\mathcal{F}}({\mathbb{N}}; \mathbb{R}^{m})$. Then, with the help of (\ref{16}), one has
\begin{eqnarray*}
  J^0(0, u) &\hspace{-0.8em}=&\hspace{-0.8em} \mathbb{E}% \sum\limits^{N-1}_{t=0}[\tilde{x}'_tQ_t\tilde{x}_t+2u'_tS_t\tilde{x}_t+u'_tR_tu_t]
%  +\tilde{x}'_{N}H\tilde{x}_{N}\nonumber\\
%  &\geq&{\color{blue}\sum\limits^{N-1}_{t=0}[\tilde{x}'_t(Q_t-S'_tR^{\dag}_tS_t)\tilde{x}_t
%  +(u_t+R^{\dag}_tS_t\tilde{x}_t)'R_t(u_t+R^{\dag}_tS_t\tilde{x}_t)]},
\sum\limits^{N-1}_{t=0}(u_t-K^*_t{x}_t)'\hat{R}_t(u_t-K^*_t{x}_t)\geq 0.
\end{eqnarray*}

(2) The proof will be shown by contradiction. Suppose that $u_t\mapsto J(x_0, u)$ is not convex, then there exists a controller $\bar{u}_t$ such that $J^0(0, \bar{u}_t)< 0$.
$\forall \lambda\in \mathbb{R}$, $J(x_0, \lambda\bar{u}_t)$ is written as

\begin{eqnarray*}
J(x_0, \lambda\bar{u}_t)&\hspace{-0.8em}=&\hspace{-0.8em}J^0(x_0, 0)+\lambda^2J^0(0, \bar{u}_t)\nonumber\\
&\hspace{-0.8em}&\hspace{-0.8em}+2\lambda \mathbb{E}\Big(\sum\limits^{N-1}_{t=0}[S_tx_t+(B_t+w_tD_t)'\lambda_t]\Big),
\end{eqnarray*}
where $x_t$ and $\lambda_t$ are the FBDEs satisfying (\ref{1}) and (\ref{4}), respectively.
Taking $\lambda \rightarrow \infty$, it yields
\begin{eqnarray*}
  V(x_0) \leq \lim_{\lambda \rightarrow \infty} J(x_0, \lambda\bar{u}_t)=-\infty,
\end{eqnarray*}
which is a contradiction.

(3) Assume $x_t$ satisfies
\begin{eqnarray}\label{17-1}
x_{t+1} &\hspace{-0.8em}=&\hspace{-0.8em} (A_t+B_tK_t)x_t+B_tv_t+[(C_t+D_tK_t)x_t\nonumber\\
&\hspace{-0.8em}&\hspace{-0.8em}+D_tv_t]w_t
\end{eqnarray}
with the initial value $x_0=0$. Combining with the Riccati equation (\ref{12-2}) and substituting $u_t=K_tx_t+v_t$ into (\ref{1})-(\ref{2}) with $x_0=0$ for $t\in \mathbb{N}$, it derives
\begin{eqnarray*}
  J^0(0, Kx+v) &\hspace{-0.8em}=&\hspace{-0.8em} \mathbb{E}\Big(\sum\limits^{N-1}_{t=0}[v'_t\hat{R}_tv_t
  +2v'_t(\hat{R}_tK_t+B'_tP_{t+1}A_t\nonumber\\
  &\hspace{-0.8em}&\hspace{-0.8em}+D'_tP_{t+1}C_t+S_t)x_t]\Big).
\end{eqnarray*}

Since $J^0(0, u)\geq \alpha E\Big(\sum\limits^{N-1}_{t=0}|u_t|^2\Big)$, we have
\begin{eqnarray}\label{17}
  0 &\leq& \alpha \mathbb{E}\Big( \sum\limits^{N-1}_{t=0} x'_tK'_tK_tx_t\Big)
  =\alpha \mathbb{E}\Big(\sum\limits^{N-1}_{t=0} |K_tx_t|^2\Big)\nonumber\\
  &\hspace{-0.8em}=&\hspace{-0.8em}\alpha \mathbb{E}\Big(\sum\limits^{N-1}_{t=0}|K_tx_t+v_t|^2\Big)
  -\alpha \mathbb{E}\Big(\sum\limits^{N-1}_{t=0}(v^2_t+2v'_tK_tx_t)\Big)\nonumber\\
  &\leq&\mathbb{E}\Big(\sum\limits^{N-1}_{t=0}\{v'_t(\hat{R}_t-\alpha I))v_t
  +2v'_t[(\hat{R}_t-\alpha I)K_t\nonumber\\
  &\hspace{-0.8em}&\hspace{-0.8em}+B'_tP_{t+1}A_t+D'_tP_{t+1}C_t+S_t]x_t\}\Big).\nonumber\\
\end{eqnarray}

For any fixed $v_0\in \mathbb{R}^m$, let $v_t=v_0$. By forward iteration for $x_0=0$ in (\ref{17-1})
and denoting $\mathcal{A}_t=A_t+B_tK_t$, $\mathcal{C}_t=C_t+D_tK_t$, $\mathcal{B}_t=B_t+w_tD_t$, $H(t,l)=\prod^{t}_{i=l}(\mathcal{A}_i+w_i\mathcal{C}_i)$
for $t\geq l$ while $H(t,l)=I$ for $t<l$, we have
\begin{eqnarray*}
 x_t=\sum\limits^{t-1}_{i=0}H(t-1,i+1)\mathcal{B}_iv_0.
\end{eqnarray*}
Substituting it into (\ref{17}), it degenerates
\begin{eqnarray*}
  0
  &\hspace{-0.8em}\leq&\hspace{-0.8em}\mathbb{E}\Big(\sum\limits^{t+h}_{s=t}\{v'_0(\hat{R}_s-\alpha I))v_0+2v'_0[(\hat{R}_s -\alpha I)K_s\nonumber\\
  &\hspace{-0.8em}&\hspace{-0.8em}+B'_sP_{s+1}A_s+D'_sP_{s+1}C_s+S_s]\nonumber\\
  &\hspace{-0.8em}&\hspace{-0.8em}\times \sum\limits^{s-1}_{i=0}H(s-1,i+1)\mathcal{B}_iv_0\} \Big) .
\end{eqnarray*}
Multiplying $\frac{1}{h}$ on both sides of above inequality and taking $h\rightarrow 0$, it concludes
\begin{eqnarray}
\mathbb{E}[v'_0(\hat{R}_s-\alpha I)v_0 ]\geq 0,
\end{eqnarray}
from the arbitrariness of $v_0$, $\hat{R}_t-\alpha I \geq 0$ is obtained.
This completes the proof.

\section{Open-Loop Solvability by Perturbation Method}
For any $\varepsilon >0$, the perturbed cost function is given by
\begin{eqnarray}\label{18}
J_{\varepsilon}(x_0, u)&\hspace{-0.8em}=&\hspace{-0.8em}J(x_0, u)+\varepsilon \mathbb{E}\Big(\sum\limits^{N-1}_{t=0}|u_t|^2\Big)\nonumber\\
&\hspace{-0.8em}=&\hspace{-0.8em} \mathbb{E}\Big(\sum\limits^{N-1}_{t=0}[x'_tQ_tx_t+2u'_tS_tx_t+u'_t(R_t+\varepsilon I)u_t\nonumber\\
&\hspace{-0.8em}&\hspace{-0.8em}+2x'_tq_t+2u'_t \rho_t]+x'_{N}Hx_{N}+2x'_{N}g\Big),
\end{eqnarray}
where $x_t$ is the state subject to (\ref{1}). The problem with the perturbed cost function is given below.
\smallskip

\emph{Problem (SLQ)$_{\varepsilon}$:}
For any initial value $x_0$, find the control $\bar{u}\in L^2_{\mathcal{F}}({\mathbb{N}}; \mathbb{R}^{m})$ satisfying
\begin{eqnarray*}
 J_{\varepsilon}(x_0, \bar{u})\leq J_{\varepsilon}(x_0, u), \quad \forall u\in L^2_{\mathcal{F}}({\mathbb{N}}; \mathbb{R}^{m}).
\end{eqnarray*}

Moreover, denote the value function of the Problem  (SLQ)$_{\varepsilon}$ as
\begin{eqnarray*}
V_{\varepsilon}(x_0)=\inf_{ u_t\in L^2_{\mathcal{F}}({\mathbb{N}}; \mathbb{R}^{m})} J_{\varepsilon}(x_0, u).
\end{eqnarray*}

\begin{theorem}\label{thm3}

Problem (SLQ)$_{\varepsilon}$ is unique optimal closed-loop solvable with $(K^{\varepsilon}, v^{\varepsilon})$ satisfying
\begin{eqnarray}
 K^{\varepsilon}_t&\hspace{-0.8em}=&\hspace{-0.8em}- (R_t+B'_tP^{\varepsilon}_{t+1}B_t+D'_tP^{\varepsilon}_{t+1}D_t+\varepsilon I)^{-1}(B'_tP^{\varepsilon}_{t+1}A_t\nonumber\\
&\hspace{-0.8em}&\hspace{-0.8em}+D'_tP^{\varepsilon}_{t+1}C_t+S_t),\label{19-1}\\
 v^{\varepsilon}_t&\hspace{-0.8em}=&\hspace{-0.8em}-[R_t+B'_tP^{\varepsilon}_{t+1}B_t+D'_tP^{\varepsilon}_{t+1}D_t]^{-1}
 \{\mathbb{E}[(B_t+w_t\nonumber\\
&\hspace{-0.8em}&\hspace{-0.8em}\times D_t)'\eta^{\varepsilon}_t|
\mathcal{F}_{t-1}]+B'_tP^{\varepsilon}_{t+1}b_t+D'_tP^{\varepsilon}_{t+1}\sigma_t+\rho_t\},\label{19-2}
\end{eqnarray}
where $\eta^{\varepsilon}$ is the adapted solution to the BSDE
\begin{eqnarray}\label{19-3}
 \eta^{\varepsilon}_{t-1}&\hspace{-0.8em}=&\hspace{-0.8em}\mathbb{E}\{[(A_t+w_tC_t)'
 +({K}^{\varepsilon}_t)'(B_t+w_tD_t)']
 \eta^{\varepsilon}_t|\mathcal{F}_{t-1}\}\nonumber\\
&\hspace{-0.8em}&\hspace{-0.8em}
 +(A_t+B_t{K}^{\varepsilon}_t)'P^{\varepsilon}_{t+1}b_t+(C_t+D_t{K}^{\varepsilon}_t)'\nonumber\\
 &\hspace{-0.8em}&\hspace{-0.8em}\times P^{\varepsilon}_{t+1}\sigma_t+q_t+({K}^{\varepsilon}_t)'\rho_t
\end{eqnarray}
with the terminal value $\eta^{\varepsilon}_{N-1}=g$.
Accordingly, for the initial value $x_0$, the open-loop optimal solution to Problem (SLQ)$_{\varepsilon}$ is given by
\begin{eqnarray}\label{20}
 u^{\varepsilon}_t &\hspace{-0.8em}=&\hspace{-0.8em} K^{\varepsilon}_tx^{\varepsilon}_t+v^{\varepsilon}_t,
\end{eqnarray}
where $x^{\varepsilon}_t$ is the solution to the closed-loop system
\begin{eqnarray}\label{21}
 x^{\varepsilon}_{t+1} &\hspace{-0.8em}=&\hspace{-0.8em} (A_t+B_tK^{\varepsilon}_t)x^{\varepsilon}_t+B_tv^{\varepsilon}_t+b_t
+[(C_t+D_tK^{\varepsilon}_t)x^{\varepsilon}_t\nonumber\\
 &\hspace{-0.8em}&\hspace{-0.8em}+D_tv^{\varepsilon}_t+\sigma_t]w_t,
\end{eqnarray}
with the initial value $ x^{\varepsilon}_0=x_0$ and $P^{\varepsilon}_{t}$ is the solution to the Riccati equation
\begin{eqnarray}\label{22}
P^{\varepsilon}_t&\hspace{-0.8em}=&\hspace{-0.8em} Q_t+A'_tP^{\varepsilon}_{t+1}A_t+C'_tP^{\varepsilon}_{t+1}C_t
-[A'_tP^{\varepsilon}_{t+1}B_t+C'_tP^{\varepsilon}_{t+1}D_t\nonumber\\
 &\hspace{-0.8em}&\hspace{-0.8em}+S'_t]
[R_t+B'_tP^{\varepsilon}_{t+1}B_t+D'_tP^{\varepsilon}_{t+1}D_t+\varepsilon I]^{-1}[B'_tP^{\varepsilon}_{t+1}A_t\nonumber\\
&\hspace{-0.8em}&\hspace{-0.8em}+D'_tP^{\varepsilon}_{t+1}C_t+S_t],
\end{eqnarray}
such that $R_t+B'_tP^{\varepsilon}_{t+1}B_t+D'_tP^{\varepsilon}_{t+1}D_t \geq \varepsilon I$.
\end{theorem}

\emph{Proof:}
Following from (\ref{18}), one has
\begin{eqnarray*}
J^0_{\varepsilon}(x_0, u)=J^0(x_0, u)+\varepsilon \mathbb{E}\Big(\sum\limits^{N-1}_{t=0}|u_t|^2\Big)
\end{eqnarray*}
which implies that
\begin{eqnarray*}
J^0_{\varepsilon}( 0, u)
\geq \varepsilon \mathbb{E}\Big(\sum\limits^{N-1}_{t=0}|u_t|^2\Big).
\end{eqnarray*}
Then, based on the proof of Theorem \ref{theorem2}, the results can be obtained.

Correspondingly, the relationship between the value function is given.
\begin{lemma}\label{lemma3}
Considering the system (\ref{1}), (\ref{2}) and (\ref{18}) with the initial value $x_0$,
there holds
\begin{eqnarray*}
  \lim_{\varepsilon \rightarrow 0} V_{\varepsilon}(x_0)=V(x_0).
\end{eqnarray*}
\end{lemma}

\emph{Proof:}
From (\ref{18}) and $\varepsilon >0$, it's easily obtained that
\begin{eqnarray*}
J_{\varepsilon}(x_0, u) \geq J(x_0, u)\geq V(x_0).
\end{eqnarray*}
The arbitrariness of $u\in L^2_{\mathcal{F}}({\mathbb{N}}; \mathbb{R}^{m})$ means that $V_{\varepsilon}(x_0)\geq V(x_0)$.

For any $\delta >0$, which is unrelated with $\varepsilon$, there exists a controller $u^{\delta}\in L^2_{\mathcal{F}}({\mathbb{N}}; \mathbb{R}^{m})$ satisfying
\begin{eqnarray*}
J(x_0, u^{\delta}) \leq V(x_0)+\delta.
\end{eqnarray*}
With the definition of (\ref{18}), it yields that
\begin{eqnarray*}
V_{\varepsilon}(x_0)\leq J(x_0, u^{\delta})+\varepsilon\sum\limits^{N-1}_{t=0}|u^{\delta}_t|^2\leq V(x_0)+\delta
+\varepsilon\sum\limits^{N-1}_{t=0}|u^{\delta}_t|^2.
\end{eqnarray*}
Taking $\varepsilon \rightarrow 0$, the above inequality degrades into
\begin{eqnarray*}
V(x_0)\leq V_{\varepsilon}(x_0)\leq V(x_0)+\delta.
\end{eqnarray*}
The arbitrariness of $\delta >0$ means that $\lim_{\varepsilon \rightarrow 0} V_{\varepsilon}(x_0)=V(x_0)$. This completes the proof.

To this end, the main results of this section is given, which reveals the relationship between the open-loop solution of Problem (SLQ)  and $\{u^{\varepsilon}_t\}$.
\begin{theorem}\label{theorem3}
For the initial value $x_0$, let $u^{\varepsilon}$ be defined in (\ref{20}), which is derived from the closed-loop optimal solution ($K^{\varepsilon}, v^{\varepsilon}$) for Problem (SLQ)$_{\varepsilon}$. The following condition are equivalent:

(1) Problem (SLQ)  is open-loop solvable with the initial value $x_0$;

(2) The set $\{u^{\varepsilon}_t\}$ is bounded in the Hilbert space $L^2_{\mathcal{F}}(\mathbb{N}; \mathbb{R}^m)$, that is to say
\begin{eqnarray*}
 \sup_{\varepsilon >0} \mathbb{E}\Big(\sum\limits^{N-1}_{t=0}|u^{\varepsilon}_t|^2 \Big)< \infty;
\end{eqnarray*}

(3) The set $\{u^{\varepsilon}_t\}$ converges strongly in $L^2_{\mathcal{F}}(\mathbb{N}; \mathbb{R}^m)$ when $\varepsilon \rightarrow 0$.

When any one of the conditions (1) to (3) is met, the set $\{u^{\varepsilon}_t\}$ is convergent strongly to an open-loop optimal control of Problem (SLQ) along with the initial value $x_0$ as  $\varepsilon \rightarrow 0$.
\end{theorem}

\emph{Proof:}
(1) $\Rightarrow$ (2). Since $u^{\varepsilon}$ is the open-loop optimal control for Problem (SLQ)$_{\varepsilon}$, combining with (\ref{18}), one has
\begin{eqnarray*}
V_{\varepsilon}(x_0)&\hspace{-0.8em}=&\hspace{-0.8em}J_{\varepsilon}(x_0, u^{\varepsilon})=J(x_0, u^{\varepsilon})+\varepsilon \mathbb{E} \Big(\sum\limits^{N-1}_{t=0}|u^{\varepsilon}_t|^2\Big)\nonumber\\
 & \geq& V(x_0)+\varepsilon \mathbb{E}\Big(\sum\limits^{N-1}_{t=0}|u^{\varepsilon}_t|^2\Big).
\end{eqnarray*}
Moreover, suppose $\bar{u}$ is the open-loop optimal solution to Problem (SLQ) with the initial value $x_0$. Based on the definition (\ref{18}), it generates that
\begin{eqnarray*}
 V_{\varepsilon}(x_0)&\leq& J_{\varepsilon}(x_0, \bar{u})=J(x_0, \bar{u})+\varepsilon \mathbb{E}\Big(\sum\limits^{N-1}_{t=0}|\bar{u}_t|^2\Big)\nonumber\\
 &\hspace{-0.8em}=&\hspace{-0.8em} V(x_0)+\varepsilon \mathbb{E}\Big(\sum\limits^{N-1}_{t=0}|\bar{u}_t|^2\Big).
\end{eqnarray*}
From the discussion above, it derives ($\varepsilon>0$)
\begin{eqnarray}\label{23}
 \mathbb{E}\Big(\sum\limits^{N-1}_{t=0}|u^{\varepsilon}_t|^2\Big) \leq
 \frac{1}{\varepsilon}(V_{\varepsilon}(x_0)-V(x_0))\leq \mathbb{E}\Big(\sum\limits^{N-1}_{t=0}|\bar{u}_t|^2\Big),
\end{eqnarray}
which implies that $\sup_{\varepsilon >0} \mathbb{E}\Big(\sum\limits^{N-1}_{t=0}|u^{\varepsilon}_t|^2\Big) < \infty$.

(2) $\Rightarrow$ (1). The boundedness of $\{u^{\varepsilon}_t\} \subseteq L^2_{\mathcal{F}}(\mathbb{N}; \mathbb{R}^m)$ means that there exists a sequence $\{\varepsilon_k\}$, $t=\{1, 2, ..., \infty\}$ and $u^*\in L^2_{\mathcal{F}}(\mathbb{N}; \mathbb{R}^m)$ satisfying $\lim_{k \rightarrow \infty}\varepsilon_k=0$ and  $\{u^{\varepsilon_k}_t\}$ is convergent weakly to a $u^*_t$, i.e.,
\begin{eqnarray*}
  \|u^*_t\|\leq \lim_{k \rightarrow \infty}\inf \|u^{\varepsilon_k}_t\|.
\end{eqnarray*}
Since $J(x_0, u)$ is convex and continuous with respect to $u$, then there holds
\begin{eqnarray*}
J(x_0, u^*)&\leq& \lim_{k \rightarrow \infty}\inf J( x_0, u^{\varepsilon_k})\nonumber\\
&\hspace{-0.8em}=&\hspace{-0.8em}\lim_{k \rightarrow \infty}\inf\Big[J_{\varepsilon_k}( x_0, u^{\varepsilon_k})-\varepsilon_k \mathbb{E}\Big(\sum\limits^{N-1}_{t=0}|u^{\varepsilon_k}_t|^2\Big)\Big]\nonumber\\
&\hspace{-0.8em}=&\hspace{-0.8em}\lim_{k \rightarrow \infty}  V_{\varepsilon_k}(x_0)=V(x_0),
\end{eqnarray*}
where the last equality is established due to Lemma \ref{lemma3}.

Hence, it's concluded that $u^*_t$ is an open-loop optimal control for the Problem (SLQ).

(3) $\Rightarrow$ (2). This part of the proof is obvious.

(2)$\Rightarrow$ (3). Suppose $\{u^{\varepsilon}_t\}$ converges strongly in $L^2_{\mathcal{F}}( \mathbb{N}; \mathbb{R}^m)$ when $\varepsilon \rightarrow 0$ and following from the previous proof, we have shown that there exists  $\{u^{\varepsilon_k}_t\}$ converges weakly to an open-loop optimal control.
In this part, firstly, we will show the set $\{u^{\varepsilon}_t\}$ is convergent weakly to a same open-loop optimal control for Problem (SLQ) with initial value $x_0$ as $\varepsilon\rightarrow 0$. That is, if there exist two sequences $\{u^{\varepsilon_{i,k}}_t\}\subseteq \{u^{\varepsilon}_t\}$ ($i=1, 2$),
which converge weakly to $u^{i, *}_t$ for $t=\{0, 1, ..., N-1\}$, then $u^{1, *}=u^{2, *}$ is an open-loop optimal
control for Problem (SLQ).

From (2) $\Rightarrow$ (1), it's obvious that $u^{1, *}$ and $u^{2, *}$ are both the open-loop optimal solution to
Problem (SLQ) with the initial value $x_0$. Since $J(x_0, u)$ is convex with respect to $u$, there holds that
\begin{eqnarray*}
J(x_0, \frac{1}{2}(u^{1, *}+u^{2, *}))  \leq  \frac{1}{2} J(x_0, u^{1, *})
+     \frac{1}{2} J( x_0, u^{2, *}) =V(x_0).
\end{eqnarray*}
That is, $\frac{1}{2}(u^{1, *}+u^{2, *})$ is also an open-loop optimal solution with
the initial value  $x_0$ to Problem (SLQ).  By using (\ref{23}), it yields
\begin{eqnarray}
 \mathbb{E}\Big(\sum\limits^{N-1}_{t=0}|u^{\varepsilon_{i,k}}_t|^2\Big) \leq
\mathbb{E}\Big[\sum\limits^{N-1}_{t=0}|\frac{1}{2}(u^{1, *}_t+u^{2, *}_t)|^2\Big], \quad i=1, 2.
\end{eqnarray}
Due to $\|u^{i, *}_t\|\leq \lim_{k \rightarrow \infty}\inf \|u^{\varepsilon_{i,k}}_t\|$, the above inequality degenerates into the following form
\begin{eqnarray}
 \mathbb{E}\Big(\sum\limits^{N-1}_{t=0}|u^{i, *}_t|^2\Big) \leq
 \mathbb{E}\Big[\sum\limits^{N-1}_{t=0}|\frac{1}{2}(u^{1, *}_t+u^{2, *}_t)|^2\Big], \quad i=1, 2,
\end{eqnarray}
which generates
\begin{eqnarray*}
2\mathbb{E}\Big[\sum\limits^{N-1}_{t=0}(u^{1, *}_t)^2+\sum\limits^{N-1}_{t=0}(u^{2, *}_t)^2\Big]\leq
\mathbb{E}\Big[\sum\limits^{N-1}_{t=0}(u^{1, *}_t+u^{2, *}_t)^2\Big],
\end{eqnarray*}
i.e., $\mathbb{E}\Big(\sum\limits^{N-1}_{t=0}(u^{1, *}_t-u^{2, *}_t)^2\Big)\leq 0$. It is established if and only if $u^{1, *}_t=u^{2, *}_t$ with $t=\{0, 1, ..., N-1\}$. That is  to say, the set $\{u^{\varepsilon}_t\}$ converge weakly to an open-loop optimal control of Problem (SLQ) along with the initial value $x_0$ as  $\varepsilon \rightarrow 0$

Next, the set $\{u^{\varepsilon}_t\}$ converges strongly in $L^2_{\mathcal{F}}(\mathbb{N}; \mathbb{R}^m)$ when $\varepsilon \rightarrow 0$ will be proved.
Suppose $u^*$ is an open-loop optimal control of Problem (SLQ), then $\{u^{\varepsilon}_t\}$ converge weakly to $u^*_t$ for $t=\{0, 1, ..., N-1\}$, i.e., $\|u^*_t\|\leq \lim_{\varepsilon \rightarrow 0}\inf \|u^{\varepsilon}_t\|$, which means that
\begin{eqnarray}\label{24}
\mathbb{E}\Big(\sum\limits^{N-1}_{t=0}|u^*_t|^2\Big)\leq \lim_{\varepsilon \rightarrow 0}\inf \mathbb{E}\Big(\sum\limits^{N-1}_{t=0}|u^{\varepsilon}_t|^2\Big).
\end{eqnarray}
Moreover, from (\ref{23}) there always holds
\begin{eqnarray}\label{25}
 \mathbb{E}\Big(\sum\limits^{N-1}_{t=0}|u^{\varepsilon}_t|^2\Big) \leq
 \mathbb{E}\Big(\sum\limits^{N-1}_{t=0}|u^*_t|^2\Big).
\end{eqnarray}
Combining (\ref{24}) with (\ref{25}), we have
\begin{eqnarray*}
\lim_{\varepsilon \rightarrow 0}\mathbb{E}\Big(  \sum\limits^{N-1}_{t=0}|u^{\varepsilon}_t|^2\Big)
=\mathbb{E}\Big(\sum\limits^{N-1}_{t=0}|u^*_t|^2\Big).
\end{eqnarray*}

Based on the discussion above, it yields
\begin{eqnarray*}
\lim_{\varepsilon \rightarrow 0}\mathbb{E}\Big(\sum\limits^{N-1}_{t=0}(u^{\varepsilon}_t-u^*_t)^2\Big)
&\hspace{-0.8em}=&\hspace{-0.8em}\lim_{\varepsilon \rightarrow 0}\mathbb{E}\Big[\sum\limits^{N-1}_{t=0}|u^{\varepsilon}_t|^2+\sum\limits^{N-1}_{t=0}|u^*_t|^2\nonumber\\
&\hspace{-0.8em}&\hspace{-0.8em}
-2\sum\limits^{N-1}_{t=0}(u^*_t)'u^{\varepsilon}_t\Big]=0.
\end{eqnarray*}
That is to say, $\{u^{\varepsilon}_t\}$ converges strongly to $u^*_t$ as $\varepsilon \rightarrow 0$.

\section{Weak Closed-Loop Solvability of Problem (LQ)}
It has been shown that if $Range(B'_tP_{t+1}A_t+D'_tP_{t+1}C_t+S_t) \nsubseteq  Range( \hat{R}_t)$, then the gain matrix $K$  in (\ref{9-1}) calculated by the Riccati equation (\ref{7}) will no-longer be the closed-loop optimal solution of Problem (SLQ). In this section, the relation with the open-loop and weak closed-loop solution to Problem (SLQ) is revealed.

Firstly, the solvability relation for Problem (SLQ) and Problem (SLQ)$^0$ is revealed.
\begin{lemma}\label{L4}
Under the condition of Assumption \ref{A1}-\ref{A2}, if Problem (SLQ) is open-loop solvable,
then Problem (SLQ)$^0$ is also open-loop solvable.
\end{lemma}

\emph{Proof:}
Following from Theorem \ref{theorem3}, if the associated set $\{u^{\varepsilon}_t\}\subseteq L^2_{\mathcal{F}}(\mathbb{N}; \mathbb{R}^m)$ in Problem (SLQ)$^0$ is bounded, then Problem (SLQ)$^0$ is open-loop solvable with the initial value $x_0$.

In Problem (SLQ)$^0$, $b_t=\sigma_t=q_t=\rho_t=g=0$, repeating the proof of Theorem \ref{theorem1}, it obtains that $\eta^{\varepsilon}$ in (\ref{8-4}) with the terminal value $\eta^{\varepsilon}_{N-1}=g=0$ is zero. Accordingly, by (\ref{13-1}), it derives $v^{\varepsilon}_t\equiv 0$. Thus, $u^{\varepsilon}_t=K^{\varepsilon}_tx^{\varepsilon}_t$ for $t=\{0, 1, ..., N-1\}$ where $x^{\varepsilon}_t$ satisfies
\begin{eqnarray}\label{25-1}
x^{\varepsilon}_{t+1} &\hspace{-0.8em}=&\hspace{-0.8em} (A_t+B_tK^{\varepsilon}_t)x^{\varepsilon}_t
+[(C_t+D_tK^{\varepsilon}_t)x^{\varepsilon}_t]w_t,
\end{eqnarray}
with the initial value $x^{\varepsilon}_0=x_0$.

In Problem (SLQ), let $u^{\varepsilon}$ be defined in (\ref{20}), where $(K^{\varepsilon}, v^{\varepsilon})$ defined in (\ref{19-1})-(\ref{19-2}) are unrelated with the initial value $x_0$. Suppose $x^{x_0, \varepsilon}_t$ and $x^{0, \varepsilon}_t$ are the solutions to (\ref{21}) with the initial value $x^{x_0,\varepsilon}_0=x_0$ and $x^{0, \varepsilon}_0=0$ on $t=\{0, 1, ..., N-1\}$, respectively. If Problem (SLQ) is open-loop solvable with the initial value $x_0$ and $0$, according to Theorem  \ref{thm3} and Theorem \ref{theorem3},
$u^{x_0, \varepsilon}=K^{\varepsilon}x^{x_0, \varepsilon}+v^{\varepsilon}$ and
$u^{0, \varepsilon}=K^{\varepsilon}x^{0, \varepsilon}+v^{\varepsilon}$ belong to $L^2_{\mathcal{F}}(\mathbb{N}; \mathbb{R}^m)$ are bounded. Since $x^{x_0, \varepsilon}_t-x^{0, \varepsilon}_t$ satisfies
\begin{eqnarray*}
x^{x_0, \varepsilon}_{t+1}-x^{0, \varepsilon}_{t+1}&\hspace{-0.8em}=&\hspace{-0.8em}(A_t+B_tK^{\varepsilon}_t)(x^{x_0, \varepsilon}_t-x^{0, \varepsilon}_t)
+[(C_t+D_tK^{\varepsilon}_t)\nonumber\\
&\hspace{-0.8em}&\hspace{-0.8em}    \times (x^{x_0, \varepsilon}_t-x^{0, \varepsilon}_t)]w_t,
\end{eqnarray*}
with the initial value $x^{x_0, \varepsilon}_0-x^{0, \varepsilon}_0=x_0$, i.e., $x^{x_0, \varepsilon}_t-x^{0, \varepsilon}_t$ is the solution to (\ref{25-1}). Thus, $u^{\varepsilon}_t$ in Problem (SLQ)$^0$ for $t=\{0, 1, ..., N-1\}$ is exactly such that
$u^{\varepsilon}_t=K^{\varepsilon}_t(x^{x_0, \varepsilon}_t-x^{0, \varepsilon}_t)=u^{x_0, \varepsilon}_t-u^{0, \varepsilon}_t$, which is bounded due to the boundedness of $\{u^{x_0, \varepsilon}_t\}$ and $\{u^{0, \varepsilon}_t\}$. Combining with the results in Theorem
\ref{theorem3}, that is to say, Problem (SLQ)$^0$ is open-loop solvable.

Next, $K^{\varepsilon}$ and $v^{\varepsilon}$ defined in (\ref{19-1})-(\ref{19-2}) will be shown that they are convergent locally in $\tilde{\mathbb{N}}$, respectively.
\begin{lemma}\label{lemma4}
Under Assumption \ref{A1}-\ref{A2}, assume Problem (SLQ)$^0$ is open-loop solvable. Then, for any $\tilde{\mathbb{N}} \subseteq \mathbb{N}$, the set $\{K^{\varepsilon}_t\}$ obtained in (\ref{19-1}) is convergent locally in $L^2(\tilde{\mathbb{N}}; \mathbb{R}^{m\times n})$, i.e.,
\begin{eqnarray}\label{28}
\lim_{\varepsilon \rightarrow 0}\sum\limits_{t\in \tilde{\mathbb{N}}}|K^{\varepsilon}_t-K^*_t|^2=0, \quad  \tilde{\mathbb{N}} \subseteq \mathbb{N},
\end{eqnarray}
where $K^*\in L^2(\tilde{\mathbb{N}}; \mathbb{R}^{m\times n})$ is deterministic locally square-integrable.
\end{lemma}

\emph{Proof:}
Denote $\bar{\mathcal{A}}_t=A_t+B_tK^{\varepsilon}_t$, $\bar{\mathcal{C}}_t=C_t+D_tK^{\varepsilon}_t$,  $H_{\varepsilon}(t,l)=\prod^{t}_{i=l}(\bar{\mathcal{A}}_i+w_i\bar{\mathcal{C}}_i)$
for $t\geq l$ while $H_{\varepsilon}(t,l)=I$ for $t<l$.
%
%$H_{\varepsilon}(t, l)=(A_t+B_tK^{\varepsilon}_t)(A_{t-1}+B_{t-1}K^{\varepsilon}_{t-1})
%\cdots(A_l+B_lK^{\varepsilon}_l)$ for $k\geq l$ while $H_{\varepsilon}(t,l)=I$ for $k< l$.
Following from (\ref{25-1}), $x^{\varepsilon}_t$ can be recalculated as
\begin{eqnarray*}
x^{\varepsilon}_t=H_{\varepsilon}(t-1, 0)x_0, \quad t\in\mathbb{N},
\end{eqnarray*}
with $x^{\varepsilon}_0=x_0$.

According to Theorem \ref{theorem3}, the open-loop solvability of Problem (SLQ)$^0$ means that
$\{u^{\varepsilon}_t=K^{\varepsilon}_tx^{\varepsilon}_t\}$ defined in (\ref{20}) is strongly convergent, that is,
\begin{eqnarray*}
  u^{\varepsilon}_t =K^{\varepsilon}_tx^{\varepsilon}_t=K^{\varepsilon}_tH_{\varepsilon}(t-1, 0)x_0, \quad t\in\mathbb{N},
\end{eqnarray*}
converges strongly in $L^2_{\mathcal{F}}(\mathbb{N}; \mathbb{R}^m)$ for any $x_0$ with $\varepsilon >0$.
Denoting $\Lambda^{\varepsilon}_t=K^{\varepsilon}_tH_{\varepsilon}(t-1, 0)$, it follows that $\Lambda^{\varepsilon}$ converges strongly in $L^2_{\mathcal{F}}(\mathbb{N}; \mathbb{R}^{m\times n})$, and the strong limits is denoted as $\Lambda^{*}$, which indicates that
\begin{eqnarray}\label{26}
\lim_{\varepsilon \rightarrow 0 }\mathbb{E}\Big(\sum\limits^{N-1}_{t=0}|\Lambda^{\varepsilon}_t-\Lambda^{*}_t|^2\Big)=0, \quad \varepsilon>0.
\end{eqnarray}

Next, the convergence of $E[H_{\varepsilon}(t, 0)]$ is proved. By the definition, $\mathbb{E}[H_{\varepsilon}(t, 0)]$ is rewritten as
\begin{eqnarray*}
\mathbb{E}[H_{\varepsilon}(t, 0)]&\hspace{-0.8em}=&\hspace{-0.8em}(A_t+B_tK^{\varepsilon}_t)\mathbb{E}[H_{\varepsilon}(t-1, 0)]\nonumber\\
&\hspace{-0.8em}=&\hspace{-0.8em}A_t\mathbb{E}[H_{\varepsilon}(t-1, 0)]+B_t \mathbb{E}[\Lambda^{\varepsilon}_t]\nonumber\\
&\hspace{-0.8em}=&\hspace{-0.8em}...=\Big(\prod\limits^{t}_{i=0}A_i\Big) \mathbb{E}[H_{\varepsilon}(-1, 0)]+\sum\limits^{t}_{i=0}\Big(\prod\limits^{t}_{l=i+1}A_l\Big) \nonumber\\
&\hspace{-0.8em}&\hspace{-0.8em}\times B_i\mathbb{E}[\Lambda^{\varepsilon}_i],
\end{eqnarray*}
with the initial value $\mathbb{E}[H_{\varepsilon}(-1, 0)]=I$. With the aid of H$\ddot{o}$lder inequality, one has
\begin{eqnarray*}
  \sum^{N-1}_{t=0}|\mathbb{E}[\Lambda^{\varepsilon}_t]-\mathbb{E}[\Lambda^{*}_t]|^2 \leq
  \mathbb{E}\Big(\sum^{N-1}_{t=0}|\Lambda^{\varepsilon}_t-\Lambda^{*}_t|^2 \Big ).
\end{eqnarray*}
Combining with (\ref{26}), it obtains
\begin{eqnarray*}
\lim_{\varepsilon \rightarrow 0}\sum^{N-1}_{t=0}|\mathbb{E}[\Lambda^{\varepsilon}_t]-\mathbb{E}[\Lambda^{*}_t]|^2=0.
\end{eqnarray*}
That is to say, $E[H_{\varepsilon}(t, 0)]$ is convergent uniformly to the solution of
\begin{eqnarray*}
\mathbb{E}[H^{*}_{\varepsilon}(t, 0)]
=A_t\mathbb{E}[H^{*}_{\varepsilon}(t-1, 0)]+B_t \mathbb{E}[\Lambda^{*}_t],
\end{eqnarray*}
with the initial value $E[H^{*}_{\varepsilon}(-1, 0)]=I$.

By using $\mathbb{E}[H^{*}_{\varepsilon}(-1, 0)]=I$, there exists a small constant $\Delta_s>0$, for $t\in [s, s+\Delta_s]$, satisfying $\mathbb{E}[H_{\varepsilon}(t-1, 0)]$ is invertible and $|\mathbb{E}[H_{\varepsilon}(t-1, 0)]|^2\geq \frac{1}{3}$.
To this end, we have
\begin{eqnarray}\label{27}
&\hspace{-0.8em}&\hspace{-0.8em}\sum\limits^{s+\Delta_s}_{t=s}|K^{\varepsilon_1}_t
-K^{\varepsilon_2}_t|^2\nonumber\\
=&\hspace{-0.8em}&\hspace{-0.8em}\sum\limits^{s+\Delta_s}_{t=s}|\mathbb{E}[\Lambda^{\varepsilon_1}_t] \mathbb{E}[H_{\varepsilon_1}(k-1, 0)]^{-1}\nonumber\\
&\hspace{-0.8em}&\hspace{-0.8em}
-\mathbb{E}[\Lambda^{\varepsilon_2}_t] \mathbb{E}[H_{\varepsilon_1}(k-1, 0)]^{-1}+\mathbb{E}[\Lambda^{\varepsilon_2}_t] \mathbb{E}[H_{\varepsilon_1}(k-1, 0)]^{-1}\nonumber\\
&\hspace{-0.8em}&\hspace{-0.8em}-\mathbb{E}[\Lambda^{\varepsilon_2}_t] \mathbb{E}[H_{\varepsilon_2}(k-1, 0)]^{-1}|^2\nonumber\\
\leq &\hspace{-0.8em}&\hspace{-0.8em}2\sum\limits^{s+\Delta_s}_{t=s}|\mathbb{E}[\Lambda^{\varepsilon_1}_t-\Lambda^{\varepsilon_2}_t]|^2
|\mathbb{E}[H_{\varepsilon_1}(k-1, 0)]^{-1}|^2\nonumber\\
&\hspace{-0.8em}&\hspace{-0.8em}
+2\sum\limits^{s+\Delta_s}_{t=s}|\mathbb{E}[\Lambda^{\varepsilon_2}_t]|^2
|\mathbb{E}[H_{\varepsilon_1}(k-1, 0)]^{-1}\nonumber\\
&\hspace{-0.8em}&\hspace{-0.8em} -\mathbb{E}[H_{\varepsilon_2}(k-1, 0)]^{-1}|^2\nonumber\\
\leq &\hspace{-0.8em}&\hspace{-0.8em} 18\sum\limits^{s+\Delta_s}_{t=s}|\mathbb{E}[\Lambda^{\varepsilon_1}_t-\Lambda^{\varepsilon_2}_t]|^2
+2\sum\limits^{s+\Delta_s}_{t=s}|\mathbb{E}[\Lambda^{\varepsilon_2}_t]|^2
\nonumber\\
&\hspace{-0.8em}&\hspace{-0.8em}\times|\mathbb{E}[H_{\varepsilon_1}(k-1, 0)]^{-1}|^2
|\mathbb{E}[H_{\varepsilon_1}(k-1, 0)]\nonumber\\
&\hspace{-0.8em}&\hspace{-0.8em}-\mathbb{E}[H_{\varepsilon_2}(k-1, 0)]|^2
|\mathbb{E}[H_{\varepsilon_2}(k-1, 0)]^{-1}|^2\nonumber\\
\leq &\hspace{-0.8em}&\hspace{-0.8em} 18\sum\limits^{s+\Delta_s}_{t=s}|\mathbb{E}[\Lambda^{\varepsilon_1}_t-\Lambda^{\varepsilon_2}_t]|^2
+162\sum\limits^{s+\Delta_s}_{t=s}|\mathbb{E}[\Lambda^{\varepsilon_2}_t]|^2\nonumber\\
&\hspace{-0.8em}&\hspace{-0.8em}\times
\sup_{t\in[s, s+\Delta_s]}|\mathbb{E}[H_{\varepsilon_1}(k-1, 0)]-\mathbb{E}[H_{\varepsilon_2}(k-1, 0)]|^2.
\end{eqnarray}
Due to (\ref{26}) and the uniformly convergence of $E[H_{\varepsilon}(t, 0)]$, by letting $\varepsilon_1, \varepsilon_2 \rightarrow 0$, (\ref{27}) is reduced into
\begin{eqnarray*}
  \lim_{\varepsilon_1, \varepsilon_2 \rightarrow 0}\sum\limits^{s+\Delta_t}_{t=s}|K^{\varepsilon_1}_t-K^{\varepsilon_2}_t|^2=0.
\end{eqnarray*}

Due to the compactness of $\tilde{\mathbb{N}}$, finite numbers $s_1$, ... $s_r$ can be selected such that $\tilde{\mathbb{N}}\subseteq \bigcup\limits^{r}_{j=1}[s_j, s_j+\Delta_{s_j}]$ and there is  $\lim_{\varepsilon_1, \varepsilon_2\rightarrow 0} \sum\limits^{s_j+\Delta_{s_j}}_{t=s_j}|K^{\varepsilon_1}_t-K^{\varepsilon_2}_t|^2=0$.

Then, when $\varepsilon_1, \varepsilon_2 \rightarrow 0$ it derives that
\begin{eqnarray*}
\sum\limits_{t\in \tilde{\mathbb{N}}}|K^{\varepsilon_1}_t-K^{\varepsilon_2}_t|^2\leq
\sum\limits^{r}_{j=1}\sum\limits^{s_j+\Delta_{s_j}}_{t=s_j}
|K^{\varepsilon_1}_t-K^{\varepsilon_2}_t|^2\longrightarrow 0.
\end{eqnarray*}
This completes the proof.

The locally convergence on $\tilde{\mathbb{N}}$ for $v^{\varepsilon}_t$ in (\ref{19-2}) is shown accordingly.
\begin{lemma}\label{lemma5}
Under Assumption \ref{A1}-\ref{A2}, if Problem (SLQ) is open-loop solvable, then, for any $\tilde{\mathbb{N}} \subseteq \mathbb{N}$, the set $\{v^{\varepsilon}_t\}$ obtained in (\ref{19-2}) is convergent locally in $L^2_{\mathcal{F}}(\tilde{\mathbb{N}}; \mathbb{R}^{m})$, i.e.,
\begin{eqnarray}\label{28-1}
\lim_{\varepsilon \rightarrow 0}\mathbb{E}\sum\limits_{t\in \tilde{\mathbb{N}}}|v^{\varepsilon}_t-v^*_t|^2=0, \quad \tilde{\mathbb{N}} \subseteq \mathbb{N},
\end{eqnarray}
where $v^*\in L^2_{\mathcal{F}}(\tilde{\mathbb{N}}; \mathbb{R}^{m})$ is deterministic locally square-integrable.
\end{lemma}

\emph{Proof:}
Following from Theorem \ref{theorem3}, the open-loop solvability of Problem (SLQ) means that
$\{u^{\varepsilon}_t=K^{\varepsilon}_tx^{\varepsilon}_t+v^{\varepsilon}_t\}$ ($t\in\mathbb{N}$) defined in (\ref{20}) converges strongly in  $L^2_{\mathcal{F}}(\mathbb{N}; \mathbb{R}^m)$, where $x^{\varepsilon}_t$ is the solution to the closed-loop system (\ref{21}), i.e.,
\begin{eqnarray}\label{28-3}
\lim_{\varepsilon_1, \varepsilon_2 \rightarrow 0} \mathbb{E}\Big(\sum\limits^{N-1}_{t=0}|u^{\varepsilon_1}_t-u^{\varepsilon_2}_t|^2\Big)=0.
\end{eqnarray}
Consider the following inequality where $x^{\varepsilon_i}_t$ satisfies (\ref{21}) with $u^{\varepsilon_i}_t$:
\begin{eqnarray}\label{28-2}
&\hspace{-0.8em}&\hspace{-0.8em}\mathbb{E} \Big(\sum\limits_{t\in \tilde{\mathbb{N}} }|K^{\varepsilon_1}_tx^{\varepsilon_1}_t
-K^{\varepsilon_2}_tx^{\varepsilon_2}_t|^2\Big)\nonumber\\
\leq &\hspace{-0.8em}&\hspace{-0.8em}2\mathbb{E} \Big(\sum\limits_{t\in \tilde{\mathbb{N}}}|K^{\varepsilon_1}_t-K^{\varepsilon}_2|^2
|x^{\varepsilon_1}_t|^2\Big)
+2\mathbb{E} \Big(\sum\limits_{t\in \tilde{\mathbb{N}}}|K^{\varepsilon_2}_t|^2
|x^{\varepsilon_1}_t\nonumber\\
 &\hspace{-0.8em}&\hspace{-0.8em}-x^{\varepsilon_2}_t|^2\Big)\nonumber\\
\leq &\hspace{-0.8em}&\hspace{-0.8em}2\sum\limits_{t\in \tilde{\mathbb{N}}}|K^{\varepsilon_1}_t-K^{\varepsilon}_2|^2
\cdot \mathbb{E} \Big(\sup_{t\in\tilde{\mathbb{N}}}|x^{\varepsilon_1}_t|^2\Big)
+2\sum\limits_{t\in \tilde{\mathbb{N}}}|K^{\varepsilon_2}_t|^2\nonumber\\
 &\hspace{-0.8em}&\hspace{-0.8em}\times \mathbb{E} \Big(\sup_{t\in\tilde{\mathbb{N}}}
|x^{\varepsilon_1}_t-x^{\varepsilon_2}_t|^2\Big).
\end{eqnarray}

Next, it shows that
\begin{eqnarray}\label{29-1}
E\Big(\sup_{t\in \mathbb{N}}|x_t|^2\Big) &\hspace{-0.8em}\leq &\hspace{-0.8em}L\cdot E \Big(|x_0|^2+ \sum^{N-1}_{t=0} |u_t|^2
+ \sum^{N-1}_{t=0}|b_t|^2\nonumber\\
 &\hspace{-0.8em}&\hspace{-0.8em}+ \sum^{N-1}_{t=0}|\sigma_t|^2\Big),
\end{eqnarray}
where $x_t$, $u_t$ are subject to (\ref{1}) and $L>0$ represents a genetic constant which can be different from line to line. By iteration, (\ref{1}) is rewritten as
\begin{eqnarray*}
x_t&\hspace{-0.8em}=&\hspace{-0.8em}\Big(\prod\limits^{t-1}_{i=0}\bar{A}_i\Big) x_0+\sum\limits^{t-1}_{i=0}\Big(\prod\limits^{t-1}_{l=i+1}\bar{A}_l\Big) \bar{B}_i[u_i+s_i],
\end{eqnarray*}
where $\bar{A}_t=A_t+w_tC_t$, $\bar{B}_t=B_t+w_tD_t$ and $s_t=b_t+w_t\sigma_t$.

Thus,  constants $L>0$ can be selected  such that
\begin{eqnarray*}
\mathbb{E}\Big(\sup_{t\in \mathbb{N}}|x_t|^2)&\hspace{-0.8em}\leq&\hspace{-0.8em} L\cdot \mathbb{E}\Big(|x_0+\sum\limits^{t-1}_{i=0}[u_i+s_i]|^2\Big)\nonumber\\
&\hspace{-0.8em}\leq&\hspace{-0.8em} L \cdot \mathbb{E}\Big((|x_0|^2+|\sum\limits^{t-1}_{i=0}u_i|^2+|\sum\limits^{t-1}_{i=0}s_i|^2)\Big)\nonumber\\
&\hspace{-0.8em}\leq&\hspace{-0.8em} L\cdot \mathbb{E}\Big((|x_0|^2+\sum\limits^{N-1}_{i=0}|u_i|^2)\nonumber\\
&\hspace{-0.8em}&\hspace{-0.8em}+ \sum^{N-1}_{t=0}|b_t|^2+ \sum^{N-1}_{t=0}|\sigma_t|^2\Big).
\end{eqnarray*}

According to (\ref{29-1}), it follows that
\begin{eqnarray}\label{29-2}
\mathbb{E}\Big(\sup_{t\in \mathbb{N}}|x^{\varepsilon_1}_t-x^{\varepsilon_2}_t|^2\Big) \leq
L \cdot \mathbb{E}\Big(\sum^{N-1}_{t=0} |u^{\varepsilon_1}_t-u^{\varepsilon_2}_t|^2\Big).
\end{eqnarray}
Based on the discussion in (\ref{28-3}), we have
\begin{eqnarray}\label{29-3}
\lim_{\varepsilon_1, \varepsilon_2\rightarrow 0}\mathbb{E}\Big(\sup_{t\in \mathbb{N}}|x^{\varepsilon_1}_t-x^{\varepsilon_2}_t|^2\Big)=0.
\end{eqnarray}
With the aid of Lemma \ref{L4} and Lemma \ref{lemma4}, reconsidering (\ref{28-2}), it degenerates into
\begin{eqnarray}\label{29-3}
\lim_{\varepsilon_1, \varepsilon_2\rightarrow 0} \mathbb{E} \Big(\sum\limits_{t\in \tilde{\mathbb{N}}}|K^{\varepsilon_1}_tx^{\varepsilon_1}_t
-K^{\varepsilon_2}_tx^{\varepsilon_2}_t|^2\Big)=0.
\end{eqnarray}

Hence, we have
\begin{eqnarray*}
&\hspace{-0.8em}&\hspace{-0.8em} \mathbb{E} \Big(\sum\limits_{t\in \tilde{\mathbb{N}}}|v^{\varepsilon_1}_t-v^{\varepsilon_2}_t|^2\Big)\nonumber\\
=&\hspace{-0.8em}&\hspace{-0.8em}\mathbb{E} \Big(\sum\limits_{t\in \tilde{\mathbb{N}}}|[u^{\varepsilon_1}_t-K^{\varepsilon_1}_tx^{\varepsilon_1}_t]
-[u^{\varepsilon_2}_t-K^{\varepsilon_2}_tx^{\varepsilon_2}_t]|^2\Big)\nonumber\\
\leq&\hspace{-0.8em}&\hspace{-0.8em}2 \mathbb{E} \Big(\sum\limits_{t\in \tilde{\mathbb{N}}}|[u^{\varepsilon_1}_t-u^{\varepsilon_2}_t|^2\Big)
+2\mathbb{E} \Big(\sum\limits_{t\in \tilde{\mathbb{N}}}|K^{\varepsilon_1}_tx^{\varepsilon_1}_t
-K^{\varepsilon_2}_tx^{\varepsilon_2}_t|^2\Big).
\end{eqnarray*}
From (\ref{28-3}) and (\ref{29-3}), it derives
\begin{eqnarray*}
\lim_{\varepsilon_1, \varepsilon_2\rightarrow 0}\mathbb{E} \Big(\sum\limits_{t\in \tilde{\mathbb{N}}}|v^{\varepsilon_1}_t-v^{\varepsilon_2}_t|^2\Big)=0.
\end{eqnarray*}
This completes the proof.

To this end, now, we are in position to propose the equivalent relationship between the open-loop and weak closed-loop solvability for Problem (SLQ).
\begin{theorem}\label{theorem5}
If Problem (SLQ) is open-loop solvable, then $(K^*, v^*)$ derived in Lemma \ref{lemma4} and  Lemma \ref{lemma5} is a weak closed-loop optimal solution for Problem (SLQ) with $t\in \tilde{\mathbb{N}}$.
Accordingly, the open-loop solvability of Problem (SLQ) is equivalent to the weak closed-loop solvability of Problem (SLQ).
\end{theorem}

\emph{Proof:}
Due to the open-loop solvability of Problem (SLQ), and following from Theorem \ref{theorem3},
$\{u^{\varepsilon}_t\}$  defined in (\ref{20}) is strongly convergent to an open-loop optimal control of Problem (SLQ), which is called as $\{u^*_t\}$. The corresponding optimal state, which is called as $x^*$ is such that
\begin{eqnarray*}
  x^*_{t+1} &\hspace{-0.8em}=&\hspace{-0.8em} A_tx^*_t+B_tu^*_t+b_t+(C_tx^*_t+D_tu^*_t+\sigma_t)w_t,
\end{eqnarray*}
with the initial value $x^*_0=x_0$.

Next, we will prove that $u^*=K^*x^*+v^*$. Then $(K^*, v^*)$ is exactly a weak closed-loop solution to  Problem (SLQ).
To justify the argument, reconsidering (\ref{29-1}), it follows that
\begin{eqnarray}\label{30}
\mathbb{E}\Big(\sup_{t\in \mathbb{N}}|x^{\varepsilon}_t-x^*_t|^2\Big) \leq
L \cdot \mathbb{E}\Big(\sum^{N-1}_{t=0} |u^{\varepsilon}_t-u^*_t|^2\Big),
\end{eqnarray}
where $x^{\varepsilon}$ is the solution to (\ref{21}). Combining (\ref{30}) with Theorem \ref{theorem3}, one has
\begin{eqnarray}\label{31}
\lim_{\varepsilon\rightarrow 0}\mathbb{E}\Big(\sup_{t\in \mathbb{N}}|x^{\varepsilon}_t-x^*_t|^2\Big)=0.
\end{eqnarray}
For $t\in \tilde{\mathbb{N}}$, it follows that
\begin{eqnarray*}
&\hspace{-0.8em}&\hspace{-0.8em}\mathbb{E}\Big(\sum_{t\in \tilde{\mathbb{N}}} |u^{\varepsilon}_t-[K^*_tx^*_t+v^*_t]|^2\Big) \nonumber\\
=&\hspace{-0.8em}&\hspace{-0.8em}
\mathbb{E}\Big(\sum_{t\in \tilde{\mathbb{N}}} |[K^{\varepsilon}_tx^{\varepsilon}_t+v^{\varepsilon}_t]-[K^*_tx^*_t+v^*_t]|^2\Big)\nonumber\\
\leq&\hspace{-0.8em}&\hspace{-0.8em} 2\mathbb{E}\Big(\sum_{t\in \tilde{\mathbb{N}}} |K^{\varepsilon}_tx^{\varepsilon}_t-K^*_tx^*_t|^2\Big)
+2\mathbb{E}\Big(\sum_{t\in \tilde{\mathbb{N}}} |v^*_t-v^{\varepsilon}_t|^2\Big)\nonumber\\
\leq&\hspace{-0.8em}&\hspace{-0.8em} 4\mathbb{E}\Big(\sum_{t\in \tilde{\mathbb{N}}} |K^{\varepsilon}_t|^2 |x^{\varepsilon}_t-x^*_t|^2\Big)
+4\mathbb{E}\Big(\sum_{t\in \tilde{\mathbb{N}}} |K^{\varepsilon}_t-K^*_t|^2 |x^{\varepsilon}_t|^2\Big)\nonumber\\
&\hspace{-0.8em}&\hspace{-0.8em}+2\mathbb{E}\Big(\sum_{t\in \tilde{\mathbb{N}}} |v^*_t-v^{\varepsilon}_t|^2\Big)\nonumber\\
\leq&\hspace{-0.8em}&\hspace{-0.8em} 4\sum_{t\in \tilde{\mathbb{N}}} |K^{\varepsilon}_t|^2 \cdot \mathbb{E}\Big(\sup_{t\in\mathbb{N}} |x^{\varepsilon}_t-x^*_t|^2\Big)
+4 \sum_{t\in \tilde{\mathbb{N}}} |K^{\varepsilon}_t-K^*_t|^2 \nonumber\\
&\hspace{-0.8em}&\hspace{-0.8em}\times \mathbb{E}\Big(\sup_{t\in\mathbb{N}} |x^{\varepsilon}_t|^2\Big)+2\mathbb{E}\Big(\sum_{t\in \tilde{\mathbb{N}}} |v^*_t-v^{\varepsilon}_t|^2\Big).
\end{eqnarray*}
By using (\ref{31}) and the results in Lemma \ref{lemma4} and  Lemma \ref{lemma5},
it has
\begin{eqnarray}
\lim_{\varepsilon\rightarrow 0}\mathbb{E}\Big(\sum_{t\in \tilde{\mathbb{N}}} |u^{\varepsilon}_t-[K^*_tx^*_t+v^*_t]|^2\Big)=0.
\end{eqnarray}

As stated in Theorem \ref{theorem3}, $u^{\varepsilon}\in L^2_{\mathcal{F}}(\mathbb{N}; \mathbb{R}^m)$ converges strongly to $u^*\in L^2_{\mathcal{F}}(\mathbb{N}; \mathbb{R}^m)$ for $t\in\mathbb{N}$, which means $u^*_t=K^*_tx^*_t+v^*_t$ with $t\in \mathbb{N}$ is established. That is to say, the open-loop solvability indicates the weak close-loop solvability. Moreover, using the definition of weak close-loop solvability, the weak close-loop solvable of Problem (SLQ) must be open-loop solvable. The proof is completed.

\section{Numerical Example}
Consider the following system
\begin{eqnarray}\label{11}
x_{t+1} &\hspace{-0.8em}=&\hspace{-0.8em} (x_t+u_t)w_t, \quad t\in\{0,1\},
\end{eqnarray}
with initial value $x_0$. The associated cost function is given by
\begin{eqnarray}\label{12}
  J(x_0, u_t) &\hspace{-0.8em}=&\hspace{-0.8em} \mathbb{E}[x^2_2-\sum\limits^{1}_{t=0}u^2_t].
\end{eqnarray}
Following from the difference Riccati equation (\ref{7})
%\begin{eqnarray*}
%% \nonumber to remove numbering (before each equation)
%P_t&\hspace{-0.8em}=&\hspace{-0.8em} Q+A'P_{t+1}A-[A'P_{t+1}B+S'][R+B'P_{t+1}B]^{\dagger}[B'P_{t+1}A+S],
%\end{eqnarray*}
with the terminal value $P_2=1$, and by simple calculation, the solution of it is exactly $P_t\equiv 1$, $t\in \{0, 1, 2\}$.

It should be noted that the problem with system (\ref{11})-(\ref{12}) is open-loop solvable but not closed-loop solvable due to the following reason:
\begin{eqnarray}
 Range( R+B'P_{t+1}B+D'P_{t+1}D]&\hspace{-0.8em}=&\hspace{-0.8em}Range( 0]=\{0\},\\
 Range(B'P_{t+1}A+D'P_{t+1}C+S]&\hspace{-0.8em}=&\hspace{-0.8em}Range(1]=\mathbb{R}.
\end{eqnarray}
That is,
\begin{eqnarray*}
Range(B'P_{t+1}A+D'P_{t+1}C+S] \nsubseteq  Range( R+B'P_{t+1}B+D'P_{t+1}D].
\end{eqnarray*}

Taking $\varepsilon>0$ and by Theorem \ref{thm3}, the solution to Riccati equation (\ref{22}) is calculated as
%
%
%Based on the method in Theorem 4.4 in [Jingrui Sun], the open-loop optimal control will be find, which means there exists weak closed-loop optimal strategy and it also will be given.
\begin{eqnarray*}
P^{\varepsilon}_t&\hspace{-0.8em}=&\hspace{-0.8em} P^{\varepsilon}_{t+1}-P^{\varepsilon}_{t+1}[-1+P^{\varepsilon}_{t+1}+\varepsilon]^{-1}
P^{\varepsilon}_{t+1}\nonumber\\
&\hspace{-0.8em}=&\hspace{-0.8em} \frac{P^{\varepsilon}_{t+1}(\varepsilon-1)}{P^{\varepsilon}_{t+1}+\varepsilon-1},
\end{eqnarray*}
with the terminal value $P^{\varepsilon}_2=1$. By backward iteration, $P^{\varepsilon}_t$ is actually calculated as
\begin{eqnarray}
P^{\varepsilon}_t=\frac{\varepsilon-1}{\varepsilon+1-k}.
\end{eqnarray}
%\begin{eqnarray*}
%% \nonumber to remove numbering (before each equation)
%P^{\varepsilon}_1=\frac{\varepsilon-1}{\varepsilon},\quad
%P^{\varepsilon}_0=\frac{\varepsilon-1}{\varepsilon+1}.
%\end{eqnarray*}

Accordingly, the closed-loop optimal gain matrix is such that
\begin{eqnarray}
K^{\varepsilon}_t&\hspace{-0.8em}=&\hspace{-0.8em}-[-1+P^{\varepsilon}_{t+1}+\varepsilon ]^{-1}P^{\varepsilon}_{t+1}\nonumber\\
&\hspace{-0.8em}=&\hspace{-0.8em}-\frac{P^{\varepsilon}_{t+1}}{P^{\varepsilon}_{t+1}+\varepsilon-1}
=-\frac{P^{\varepsilon}_t}{\varepsilon-1}=-\frac{1}{\varepsilon+1-k}.
\end{eqnarray}
%i.e.,
%\begin{eqnarray*}
%% \nonumber to remove numbering (before each equation)
%K^{\varepsilon}_0 &\hspace{-0.8em}=&\hspace{-0.8em} -\frac{P^{\varepsilon}_0}{\varepsilon-1}
%=-\frac{\varepsilon-1}{\varepsilon+1}\frac{1}{\varepsilon-1}=-\frac{1}{\varepsilon+1},\\
%K^{\varepsilon}_1 &\hspace{-0.8em}=&\hspace{-0.8em} -\frac{P^{\varepsilon}_1}{\varepsilon-1}
%=-\frac{\varepsilon-1}{\varepsilon}\frac{1}{\varepsilon-1}=-\frac{1}{\varepsilon}.
%\end{eqnarray*}

Hence, the corresponding closed-loop state equation is written as
\begin{eqnarray}
x^{\varepsilon}_{t+1}=(1+K^{\varepsilon}_t)x^{\varepsilon}_tw_t, \quad x^{\varepsilon}_0=x_0.
\end{eqnarray}
The control is given by
\begin{eqnarray}
u^{\varepsilon}_t=K^{\varepsilon}_tx^{\varepsilon}_t
=\prod\limits^{t-1}_{i=0}[(1+K^{\varepsilon}_i)w_i]x^{\varepsilon}_0,
\end{eqnarray}
that is to say,
\begin{eqnarray}
u^{\varepsilon}_0 &\hspace{-0.8em}=&\hspace{-0.8em} K^{\varepsilon}_0 x^{\varepsilon}_0 =-\frac{1}{\varepsilon+1}x_0,\\
u^{\varepsilon}_1 &\hspace{-0.8em}=&\hspace{-0.8em} K^{\varepsilon}_1 x^{\varepsilon}_1
=K^{\varepsilon}_1(1+K^{\varepsilon}_0)x^{\varepsilon}_0w_0\nonumber\\
&\hspace{-0.8em}=&\hspace{-0.8em}-\frac{1}{\varepsilon} (1-\frac{1}{\varepsilon+1})x^{\varepsilon}_0w_0=-\frac{1}{\varepsilon+1}x_0w_0.
\end{eqnarray}

Since
\begin{eqnarray}
   \mathbb{E}[\sum\limits^{1}_{t=0}|u^{\varepsilon}_t|^2] &\hspace{-0.8em}=&\hspace{-0.8em} \frac{2x^2_0}{(\varepsilon+1)^2}\leq  2x^2_0,
\end{eqnarray}
that is to say, $u^{\varepsilon}_t$ is bounded. From the results in Theorem \ref{theorem3}, the open-loop optimal solution is given as
\begin{eqnarray}
   u^{*}_0 &\hspace{-0.8em}=&\hspace{-0.8em} \lim_{\varepsilon\rightarrow 0}u^{\varepsilon}_0=-x_0,\\
  u^{*}_1 &\hspace{-0.8em}=&\hspace{-0.8em} \lim_{\varepsilon\rightarrow 0}u^{\varepsilon}_1=-x_0w_0.
\end{eqnarray}

Finally, based on the discussion in Theorem \ref{theorem5}, the weak closed-loop optimal strategy $K^{*}_t$ is given by
\begin{eqnarray}
K^{*}_t &\hspace{-0.8em}=&\hspace{-0.8em} \lim_{\varepsilon\rightarrow 0}K^{\varepsilon}_t
=\lim_{\varepsilon\rightarrow 0}(\frac{1}{\varepsilon+1-k})
=-\frac{1}{1-k}.
\end{eqnarray}

\section{Conclusion}

The open-loop, closed-loop, and weak closed-loop solvability for the stochastic  discrete-time system have been discussed in this paper.
We have shown the closed-loop solvability is equivalent to the regular solvability of the generalized Riccati equation, while it also means that it's the open-loop solvable.
However, there is an example illustrates that when the stochastic LQ problem is open-loop solvable may not be closed-loop solvable.
For the stochastic LQ problem, which is merely open-loop solvable, we have found a weak closed-loop solution whose outcome is the open-loop solution to the stochastic LQ problem which implies that there is a linear feedback form of the state for the open-loop solution.
Moreover, an equivalent relationship between the open-loop and the weak closed-loop solutions to the stochastic LQ problem has been established.

\vfill

\end{document}